\documentclass[12pt]{article}
\usepackage{amsfonts,mathrsfs}
\def\d{{\rm d}}
\def\mi{{\rm i}}
\def\eps{\varepsilon}
\def\g{\mathop{\gamma}\nolimits}
\def\G{\mathop{\Gamma}\nolimits}
\def\D{\mathop{\Delta}\nolimits}
\def\Re{\mathop{\rm Re\,}\nolimits}

\def\e{\mathop{\rm e}\nolimits}
\def\Z{\mathop{\mathcal Z}\nolimits}
\def\Res{\mathop{\rm Res}\nolimits}
\def\hf{{\textstyle{1 \over 2}}}
\def\qt{{\textstyle{1 \over 4}}}
\def\fq{{\textstyle{5 \over 4}}}

\def\defi{\stackrel{\rm def}{=}}
\def\si{\!\!\! &}
\def\sf{& \!\!\!}

\title{Zeta functions for the Riemann zeros}
\author{{\bf Andr\'e Voros} \\
\\
CEA, Service de Physique Th\'eorique de Saclay\\
CNRS URA 2306\\
F-91191 Gif-sur-Yvette CEDEX (France)\\
{ E-mail : {\tt voros@spht.saclay.cea.fr}}\\
\\
and\\
\\
Institut de Math\'ematiques de Jussieu--Chevaleret\\
CNRS UMR 7586\\
Universit\'e Paris 7\\
2 place Jussieu,
F-75251 Paris CEDEX 05 (France)}

\begin{document}
\maketitle

{\abstract
A family of Zeta functions built as Dirichlet series over the Riemann zeros 
are shown to have meromorphic extensions in the whole complex plane, 
for which numerous analytical features 
(the polar structure, plus countably many special values) 
are explicitly displayed.
}
\bigskip

\centerline{\bf Fonctions Z\^eta pour les z\'eros de Riemann}

{\abstract
Certaines fonctions Z\^eta d\'efinies sur les z\'eros de Riemann, 
par une famille de s\'eries de Dirichlet, se prolongent \`a tout
le plan complexe en des fonctions m\'eromorphes, dont de nombreuses
caract\'eristiques peuvent \^etre explicit\'ees (la structure polaire,
mais aussi une infinit\'e de valeurs sp\'eciales).
}
\bigskip

\medskip


\bigskip

\section{Introduction}

This work proposes to investigate certain meromorphic functions 
defined by Dirichlet series over the nontrivial zeros $\{\rho\}$ 
of the Riemann zeta function $\zeta(s)$,
and to thoroughly compile their explicit analytical features.
If the Riemann zeros are listed in pairs as usual,
\begin{equation}
\label{ZER}
\{\rho = \hf \pm \mi \tau_k \}_{k=1,2,\ldots}, \quad
\{ \Re \, \tau_k \} \mbox{ positive and non-decreasing,}
\end{equation}
then the Dirichlet series to be mainly studied read as
\begin{equation}
\label{ZDef}
\Z(\sigma,v) \defi
\sum_{k=1}^\infty ({\tau_k}^2+v)^{-\sigma}, \qquad \Re \sigma > \hf,
\quad v> - {\tau_1}^2, 
\end{equation}
extended to meromorphic functions of $\sigma \in \mathbb C$, 
and parametrized by $v$
--- with emphasis on two cases, $v=\qt$ and especially $v=0$.
Their analysis will simultaneously yield some results for the variant series
\begin{equation}
\label{Hz}
{\mathfrak Z}(\sigma,a) \defi \sum_{k=1}^\infty ({\tau_k}+a)^{-2\sigma}, 
\qquad \Re \sigma > \hf \quad (\mbox{and, e.g., }|\arg a | \le \pi/2).
\end{equation}

Those Zeta functions are ``secondary": 
arising from the nontrivial zeros of a classic zeta function (here, $\zeta(s)$);
and ``generalized": they admit an auxiliary shift parameter
just like the Hurwitz zeta function 
($\zeta(s,a) \defi \sum_0^\infty (n+a)^{-s}$).
Such ``$\zeta$-Zeta" functions have occasionally appeared in the literature, 
but mostly through particular cases or under very specific aspects.
On the other hand, their {\sl abundance of general explicit properties\/}
seems to have been largely ignored,
although it can be revealed by quite {\sl elementary\/} means (as we will do).
And with regard to the Riemann zeros, reputed to be highly elusive quantities, 
those properties constitute additional explicit information: this is enough 
to motivate a more comprehensive treatment (and bibliography) of the subject.
The present work just aims to do that, 
in a self-contained and very concrete way, as a kind of
``All you ever wanted to know about $\zeta$-Zeta functions\ldots" handbook
(without prejudice to the usefulness of any single result by itself).

\medskip

We begin by developing the background and our motivations in greater detail.

First, if a {\sl Selberg zeta function\/} is used in place of $\zeta(s)$ 
from the start
(assuming the simplest setting of a compact hyperbolic surface $X$ here),
then the $\{ {\tau_k}^2+ \qt \}$ become the eigenvalues 
of the (positive) Laplacian on $X$, 
and the Zeta functions (\ref{ZDef}) turn into 
{\sl bona fide\/} spectral (Minakshisundaram--Pleijel) zeta functions,
for which numerous explicit results have indeed been displayed
(with the help of Selberg trace formulae:
cf. \cite{R} for $v \ge \qt$, \cite{S} for $v=\qt$, \cite{CV,Vz} for $v=0$).

Some transposition of those results to the Riemann case can then be expected, 
in view of the formal analogies between 
the trace formulae for Selberg zeta functions on the one hand,
and the Weil ``explicit formula" for $\zeta(s)$ on the other hand \cite{H}.
Indeed, a few symmetric functions over the Riemann zeros 
that resemble spectral functions have been well described, mainly 
$V(t) \defi \sum_\rho \e^{\rho t}$ (\cite{CR,G2}, \cite{JL2} chap.II).
Zeta functions like (\ref{ZDef}) have also been considered,
but almost solely to establish their meromorphic continuation 
to the whole $\sigma$-plane --- apart from the earliest occurrence we found:
a mention by Guinand \cite{G1} (see also \cite{CHA}) of the series 
$\sum_k \tau_k^{-s}$ ($\equiv \Z(s/2,v=0)$) on one side of a functional relation
(eq.(\ref{GFE}) below) arising as an instance 
of a generalized Poisson summation formula.
Later, Delsarte introduced that function again (as $\phi(s)$ in \cite{Dl})
to describe its poles qualitatively,
displaying (only) its principal polar part at $s=1$, as $(2\pi)^{-1}/(s-1)^2$;
Kurokawa \cite{K} made the same study at $v=\qt$, 
not only for $\zeta(s)$ but also for Dedekind zeta functions
and Selberg zeta functions for PSL$_2(\mathbb Z)$ [or congruence subgroups]
(then, Zeta functions like (\ref{ZDef}) occur 
within the {\sl parabolic\/} components);
and Matiyasevich \cite{MA} discussed the special values 
$\theta_n \equiv 2 \Z(n,\qt) \ (n\in {\mathbb N}^\ast)$.
Extensions in the style of the Lerch zeta function have also been studied
(\cite{F}, \cite{JL2} chap.VI).

Independently, Deninger \cite{Dn} and Schr\"oter--Soul\'e \cite{SS}
considered a {\sl different\/} Hurwitz-type family
(we keep their factor $(2\pi)^s$ just to avoid multiple notations),
\begin{equation}
\label{ZZet}
\xi (s,x) \defi (2\pi)^s \sum_\rho (x-\rho)^{-s} \qquad (\Re s >1),
\end{equation}
mainly to evaluate $\partial_s \xi (s,x)_{s=0}$ (as eq.(\ref{DZet1}) below);
earlier, Matsuoka \cite{Ms2}, then Lehmer \cite{L} had focused upon the sums 
\begin{equation}
\label{RN}
{\mathscr Z}_n \defi \sum_\rho \rho^{-n}, \quad n \in {\mathbb N}^\ast
\quad [\equiv (2\pi)^{-n} \xi(n,1), \mbox{ for } n \ne 1].
\end{equation}
Here we easily recover $\xi (s,x)$ from 
the {\sl other\/} Zeta function (\ref{Hz}) (but not the reverse), as
\begin{equation}
\label{Z2X}
\xi (s,\hf+y) \equiv (2\pi)^s 
\Bigl[ \e^{\mi\pi s/2} {\mathfrak Z} \bigl( \hf s, \mi y \bigr) +
     \e^{-\mi\pi s/2} {\mathfrak Z} \bigl( \hf s, -\mi y \bigr) \Bigr] .
\end{equation}

\bigskip

{\sl The present work proposes a broader, and unified, description 
for all those $\zeta$-Zeta functions, with a wealth of explicit results 
comparable to usual spectral zeta functions\/} \cite{Vz}.

Tools for the purpose could also be borrowed from spectral theory 
(trace formulae, etc.), but the objects under scrutiny are more singular here
(the Zeta functions for the Riemann zeros manifest double poles, 
vs simple poles in the Selberg case);
this then imposes so many adaptations upon the classic procedures
that a self-contained treatment of the Riemann case alone is actually simpler.
Even then, we cannot get maximally explicit outputs for all cases at once
(e.g., Weil's ``explicit formula" {\sl diverges\/} for $v \le \qt$),
and our analysis has to develop gradually.

So, we begin (Sec.2) by setting up a minimal abstract framework,
sufficient to handle $\Z(\sigma,v)$
(with a permanent distinction between properties in the half-planes
$\{ \Re \sigma < 1 \}$ and $\{ \Re \sigma > \hf \}$ respectively).
We next obtain a first batch of explicit results for the case $v=\qt$ (Sec.3),
then for general values of $v$ (Sec.4).
Specializing to the case $v=0$ in Sec.5, 
we reach an almost explicit meromorphic continuation formula for $\Z(\sigma,0)$
into the half-plane $\{ \Re \sigma < \hf \}$, 
which immediately implies many more properties of that function,
and is generalizable to $L$-series and other number-theoretical zeta functions.
In Sec.6 we exploit the latter results to sharpen the descriptions of 
both Hurwitz-type functions $\Z(\sigma,v)$ and ${\mathfrak Z}(\sigma,a)$.
Sec.7 provides a summary of the results; 
essentially, a Table collects the main formulae for $v=0$ and $\qt$,
also referring to the main text like an index.
(Which text can in some sense be viewed, and simply used, 
as a set of ``notes" for this Table~!)
Finally, Appendices A and B treat some subsidiary issues: 
a meromorphic continuation method for the Mellin transforms of Sec.2,
and certain numerical aspects.

\bigskip

For convenience, we recall the needed elementary results and notations 
\cite{AS,B,CHO}:
\begin{eqnarray}
\label{NOT}
B_n \si:\sf \mbox{Bernoulli numbers } (B_0=1,\ B_1=-\hf,
\ B_2=\textstyle{1 \over 6},\ldots; 
\ B_{2m+1} = 0 \mbox{ for }m=1,2,\ldots) \nonumber \\
(B_n(x) \si:\sf \mbox{Bernoulli polynomials}) \nonumber \\
\g \si=\sf \mbox{Euler's constant; \ Stieltjes constants : } 
\g_n \defi \!\! \lim_{M \to \infty} \biggl\{ 
\sum_{m=1}^M \! {(\log m)^n \over m}-{(\log M)^{n+1} \over n+1} \biggl\}
\nonumber \\
(\g_0 \si=\sf \g \approx 0.5772156649, \quad 
\g_1 \approx -0.0728158455, \quad \g_2 \approx -0.0096903632, \ldots) \\
E_n  \si:\sf \mbox{Euler numbers } (E_0=1,\ E_2=-1,\ E_4=5,\ldots;
\ E_{2m+1} = 0 \mbox{ for }m=0,1,\ldots). \nonumber 
\end{eqnarray}
[$E_n$ is also a standard  symbol for quantized energy levels,
a concept often invoked purely rhetorically about the Riemann zeros;
let us then insist that our (present) work is totally decoupled 
from such interpretations, 
whereas it sees the Euler numbers as truly essential~!]

Concerning the Riemann zeta function 
$\zeta(s) \defi \sum\limits_{k=1}^\infty k^{-s}$ \cite{T,Da,E},
we will need its special values,
\begin{equation}
\label{SVR}
\zeta(-n)= (-1)^n {B_{n+1} \over n+1} \quad (n=0,1,\ldots); \qquad 
\zeta(2m) = {(2 \pi)^{2m} \over 2(2m)!} |B_{2m}| \quad (m=1,2,\ldots);
\end{equation}
\begin{equation}
\label{Z0}
\zeta(0)=-\hf \qquad \mbox{and} \qquad \zeta'(0)=-\hf\log 2\pi,
\end{equation}
and its functional equation in the form 
\begin{equation}
\label{XI}
\Xi(s) \defi {\zeta(s) \over F(s)} \equiv {\zeta(1-s) \over F(1-s)}, \qquad
\qquad F(s) \defi {\pi^{s/2} \over s(s-1) \G(s/2)},
\end{equation}
where $\Xi(s)$ is an {\sl entire\/} function,
which is {\sl even\/} under the symmetry $s \leftrightarrow (1-s)$
(this expresses the functional equation), and only keeps the nontrivial zeros
of $\zeta(s)$.

In order to analyze $\zeta(s)$, 
we will be {\sl forced\/} to invoke a particular Dirichlet $L$-series as well
(associated with the Dirichlet character $\chi_4$ \cite{AS,B,Da})
\begin{equation}
\label{LDef}
\beta(s) \defi \sum_{k=1}^\infty (-1)^{k+1} (2k-1)^{-s} \quad [=L(\chi_4,s)];
\end{equation}
$\beta(s)$ extends to an entire function, having the special values 
\begin{equation}
\label{SVL}
\beta(-m)={E_m \over 2}, \qquad 
\beta(2m+1)= { (\pi/2)^{2m+1} \over 2 (2m)!} |E_{2m}| \qquad (m=0,1, \ldots);
\end{equation}
\begin{equation}
\label{L0}
\beta(0)= \hf \qquad \mbox{and} \qquad
\beta'(0)=-\textstyle{3 \over 2} \log 2 -\log \pi + 2 \log \G({1 \over 4}),
\end{equation}
and a functional equation expressible as
\begin{equation}
\label{xI}
\Xi_{\chi_4} (s) \defi 
\Bigl( {4 \over \pi} \Bigr)^{(s+1)/2} \G \Bigl( {s+1 \over 2} \Bigr) \beta(s)
\equiv \Xi_{\chi_4} (1-s)
\end{equation}
where the function $\Xi_{\chi_4}$ is entire 
and only keeps the nontrivial zeros of $\beta(s)$.

\section{General Delta and Zeta functions of order $<1$}

\subsection{Admissible sequences and Delta functions}

Throughout this work, a numerical sequence $\{x_k\}$
(systematically labeled by {\sl positive\/} integers $k=1,2,\ldots$) 
will be called {\sl admissible of order\/} $\mu_0<1$ if:

(i) $0 < x_1 \le x_2 \le \cdots,\quad  x_k \uparrow +\infty$ ;

\noindent (or: complex $x_k \to \infty$ with $|\arg x_k|$ sufficiently bounded
\cite{QHS,I}, to provide an unconditionally valid framework ``in case of need");

(ii) $\sum\limits_k |x_k|^{-1} < \infty$, 
making the following Weierstrass product {\sl converge\/},
\begin{equation}
\label{DDef}
\D(z \mid \{x_k\}) \defi \prod_k (1+ z /x_k) 
\qquad (\forall z \in \mathbb C) ;
\end{equation}
it then defines an {\sl entire\/} ``Delta" function $\D(z)$
(we omit the argument $\{x_k\}$ except when an ambiguity may result);

(iii) for $z \to\infty$, $\log \D(z)=o(|z|^{\mu_0+\delta}) \ \forall \delta>0 $,
and it admits a complete {\sl uniform asymptotic expansion\/}
in some sector $\{ |\arg z| < \theta \}$, of a form 
governed by some strictly decreasing sequence of real exponents $\{\mu_n\}$, as
\begin{eqnarray}
\label{DAs}
\log \D(z) \sim \sum_{n=0}^\infty 
(\tilde a_{\mu_n} \log z + a_{\mu_n})  z^{\mu_n} \qquad (z \to \infty) \\
\mbox{with } \mu_0 > \mu_1 > \cdots,\quad \mu_n  \downarrow -\infty, \quad 
\mbox{and }  0< \mu_0 < 1 \nonumber
\end{eqnarray}
(``generalized Stirling expansion", 
by extension from the case $x_k=k$ \cite{JL1});
such a uniform expansion is repeatedly differentiable in $z$.

Then, the Dirichlet series
\begin{equation}
\label{Zdef}
Z(\sigma \mid \{x_k\}) \defi \sum_{k=1}^\infty x_k^{-\sigma} \qquad 
(\mbox{convergent for } \Re \sigma > \mu_0)
\end{equation}
defines the {\sl Zeta function of\/} $\{x_k\}$, holomorphic in that half-plane.
The point $\sigma=1$ has to lie in the latter by assumption (ii), 
which imposes $\mu_0 < 1$; then $(\log \D)$ can moreover be expanded
term by term in eq.(\ref{DDef}) to yield the Taylor series
\begin{equation}
\label{DTay}
\log \D(z) \equiv 
\sum_{m=1}^\infty (-1)^{m-1} {Z(m) \over m} z^m
\qquad (\mbox{convergent for } |z|<x_1).
\end{equation}

Motivations: the idea here is to assume certain properties 
for an entire function $\D(z)$ of order $\mu_0<1$,
so as to generate a function $Z(\sigma)$ meromorphic in all of $\mathbb C$
with poles at the $\mu_n$, their maximum order $r$ being 2 here
(as dictated by the specific form of eq.(\ref{DAs})).
Such $\D(z)$ are very special instances of zeta-regularized infinite products,
for which much more general frameworks exist (e.g., \cite{JL1,JL2,I}).
However, singularities and essential complications definitely
increase each time $r$ or the integer part $[\mu_0]$ can become larger
(chiefly, the formalism leaves $([\mu_0]+1)$ meaningful integration constants
undetermined).
Efficiency then commands to minimize both latter parameters 
(subject to $r \ge 1$ and $\mu_0>0$);
specially, $\mu_0<1$ is simpler to handle than any $\mu_0 \ge 1$.
In this respect, spectral zeta functions frequently have $\mu_0 \ge 1$
(e.g., for Laplacians on compact Riemannian manifolds, 
$\mu_0 = \hf \times$[dimension]) but only simple poles,
hence the simplifying assumption $r=1$ is suitable for them \cite{V,QHS};
by contrast, the present functions $\Z(\sigma,v)$ 
will accept the lower value $\mu_0=\hf$ but will need $r=2$.
Another difference is that for eigenvalue spectra,
a ``partition function" $\sum_k \e^{-tx_k}$ is a natural starting point; 
in the Riemann case, that type of function ($V(t)=\sum_\rho \e^{\rho t}$) 
exhibits a more remote and contrived structure \cite{CR}, 
while Delta-type functions are easily accessible
(by simple alterations of eq.(\ref{XI})).
Thus, Riemann zeros and standard eigenvalue spectra 
have several {\sl mutually singular\/} features, 
making their unified handling rather cumbersome.

\subsection{Meromorphic continuation of Zeta functions}

The bounds $\log \D(z)=O(z^1)$ for $z \to 0$, and $O(z^{\mu_0})$ for $z\to+\infty$,
with $\mu_0<1$, imply that the following Mellin transform of $\log \D(z)$,
\begin{equation}
\label{IMel}
I(\sigma) \defi \int_0^{+\infty} \log \D(z) z^{-\sigma-1} \d z 
\qquad \qquad (\mu_0 < \Re \sigma <1)
\end{equation}
converges to a holomorphic function of $\sigma$ in the stated vertical strip.

Then, by standard arguments (see App.A, and \cite{J,V,C}),
$I(\sigma)$ extends to a meromorphic function on either side of that strip:

\noindent $a)$ by virtue of the expansion (\ref{DAs}) for $z \to \infty$, 
$I(\sigma)$ extends to all $\Re \sigma<1$, with 
\begin{equation}
\label{IP2}
\mbox{(at most) double poles at } \sigma=\mu_n, \quad 
\mbox{polar parts} = 
{\tilde a_{\mu_n} \over (\sigma-\mu_n)^2} + {a_{\mu_n} \over (\sigma-\mu_n)};
\end{equation}

\noindent $b)$ just because $\log \D(z)$ has a Taylor series at $z=0$,
$I(\sigma)$ extends to all $\Re \sigma>\mu_0$, with
\begin{equation}
\label{IP1}
\mbox{(at most) simple poles at } \sigma =m \in \mathbb N^\ast, \quad
\mbox{of residues } -(\log \D)^{(m)}(0)/m!\ .
\end{equation}

Now, for $\mu_0 < \Re \sigma <1$, the integral (\ref{IMel}) 
can also be done term by term after inserting the expansion (\ref{DDef}), giving
\begin{equation}
\label{ZMel}
{\sigma \sin \pi \sigma\over \pi} I(\sigma) \equiv 
\sum_{k=1}^\infty x_k^{-\sigma} \equiv Z(\sigma).
\end{equation}
The meromorphic extension of $I(\sigma)$ then entails that of $Z(\sigma)$ 
to the whole $\sigma$-plane as well, 
i.e., to the half-planes $\{ \Re \sigma <1 \}$ by $a)$,
resp. $\{ \Re \sigma > \mu_0 \}$ by $b)$ independently, so that:

\noindent $a')$ any {\sl non-integer\/} pole $(\mu_n\notin \mathbb Z)$ 
of $I(\sigma)$ generates at most a double pole for $Z(\sigma)$, with
\begin{equation}
\label{ZP2}
Z(\mu_n+ \eps) = 
\Bigl[ {\mu_n \sin \pi \mu_n \over \pi} \tilde a_{\mu_n} \Bigr] {1 \over \eps^2}
+ \Bigl[ {\mu_n \sin \pi \mu_n \over \pi} a_{\mu_n} + 
\Bigl( {\sin \pi \mu_n \over \pi} + \mu_n \cos \pi \mu_n \Bigr) \tilde a_{\mu_n} \Bigr]
{1 \over \eps} \ [+ \mbox{regular part}];
\end{equation}

\noindent $a'')$ any {\sl negative integer\/} pole $(\mu_n=-m)$ of $I(\sigma)$ 
generates at most a simple pole for $Z(\sigma)$,
through partial cancellation with the zeros of $(\sin \pi \sigma)$, with
\begin{eqnarray}
\label{ZP1}
(m \in \mathbb N^\ast): \qquad Z(-m+\eps) \si=\sf 
(-1)^m \Bigl[ -{m \tilde a_{-m} \over \eps} + (\tilde a_{-m} - m a_{-m}) \Bigr]+O(\eps)
\nonumber\\
\label{ZRFP}
i.e., \mbox{ residue:} \quad \Res_{\sigma=-m} Z(\sigma) 
\si=\sf (-1)^{m+1} m \tilde a_{-m}\ ;\\
\label{ZT}
\mbox{ finite part:} \quad {\rm FP}_{\sigma=-m} Z(\sigma) 
\si=\sf (-1)^m(\tilde a_{-m} - m a_{-m})
\end{eqnarray}
(we call eqs.(\ref{ZT}) ``trace identities" by extension from spectral theory,
specially when $\tilde a_{-m}=0$, in which case explicit finite values 
for the $Z(-m)$ result);

\noindent $a''')$ any pole at $\sigma=0$ of $I(\sigma)$ gets fully cancelled 
by the double zero of $(\sigma \sin \pi \sigma)$, giving
\begin{equation}
\label{ZP0}
Z(\eps) = \tilde a_0 + a_0 \eps + O(\eps^2) \qquad \Longrightarrow \qquad
Z(0)=\tilde a_0, \quad Z'(0)=a_0.
\end{equation}

\noindent $b')$ each {\sl positive integer\/} pole of $I(\sigma)$
($\sigma=m \in \mathbb N^\ast$, from eq.(\ref{IP1}))
gets cancelled by a zero of $(\sin \pi \sigma)$, giving
\begin{equation}
\label{ZTay}
Z(m)=(-1)^{m-1}(\log \D)^{(m)}(0)/(m-1)!
\end{equation}
but this last output just duplicates with the previous Taylor formula
(\ref{DTay}).

Ultimately, all the poles of $Z(\sigma)$ 
lie in a {\sl single\/} decreasing sequence $\{ \mu_n \}_{n \in \mathbb N}$, 
and have maximum order $r=2$ under the specific assumption (\ref{DAs}).

\section{Delta function based at $s=0$, 
and the Zeta function $Z(\sigma) \defi \Z(\sigma,v=\qt)$}

We now apply the previous framework to the Riemann zeros 
upon just slight changes with regard to the usual factorization
of $\zeta(s)$ around $s=0$.

\subsection{Basic facts and notations}

The most convenient starting point is the entire function (\ref{XI}),
which keeps precisely the non-trivial zeros $\{\rho\}$ of $\zeta(s)$, 
and has the familiar Hadamard product representation \cite{E}
\begin{equation}
\label{CHAD}
\Xi(s) = \e^{Bs} \prod_\rho (1-s /\rho) \e^{s /\rho}, 
\qquad B \defi [\log \Xi]'(0) = \log 2\sqrt\pi\ -1 -\hf\g .
\end{equation}

Zeros are henceforth grouped in pairs as 
$\{\rho = \hf \pm \mi \tau_k\}_{k=1,2,\ldots}$ according to eq.(\ref{ZER}).
Their corresponding counting function,
${\mathcal N}(T) \defi {\rm card} \{\tau_k \mid \Re \tau_k < T \}$,
follows a well-known estimate $\overline {\mathcal N}(T)$, as \cite{T,E}
\begin{equation}
\label{GE}
{\mathcal N}(T) \sim \overline {\mathcal N}(T) \defi
{T \over 2 \pi}\Bigl[ \log {T \over 2 \pi} -1 \Bigr] \qquad
(T \to +\infty).
\end{equation}

Thus, the sequence of zeros itself could be admissible of order 1 at best;
fortunately, a transformed sequence $\{x_k\}$ and its Zeta function
will immediately arise:
\begin{equation}
\label{Z1}
x_k \defi \qt +{\tau_k}^2 , \qquad {\rm and} \qquad
Z(\sigma) \defi Z(\sigma \mid \{x_k\}) = \sum_k {x_k}^{-\sigma}
\end{equation}
(the latter series will converge for $ \Re \sigma > \hf $, 
again by the estimate (\ref{GE})).

Indeed, once the zeros have been reordered in pairs, it first follows that
$ {\mathscr Z}_1 \defi \sum\limits_\rho \rho^{-1} = 
\sum\limits_k {x_k}^{-1} = Z(1) $
(convergent sums), and then,
$ \Xi(s) \equiv \e^{[B+ Z(1)]s} \prod\limits_k [1+ s(s-1)/x_k]$;  
the parity property $\Xi(s)=\Xi(1-s)$ thereupon imposes $Z(1)=-B$, hence
\begin{equation}
\label{Z11}
Z(1) = {\mathscr Z}_1 = -[\log \Xi]'(0) = -\log 2\sqrt\pi\ +1 +\hf \g ,
\end{equation}
a classic result (\cite{Da}, ch.12; \cite{E}, Sec.3.8). 
All in all, the product formulae and functional equation boil down to 
\begin{equation}
\label{PROD}
{\zeta(s) \over F(s)} \equiv \D(\lambda \mid \{x_k\}) 
\equiv {\zeta(1-s) \over F(1-s)} 
\qquad \qquad \Bigl[ F(s) \defi {\pi^{s/2} \over s(s-1) \G(s/2)} \Bigr] ,
\end{equation}
\begin{equation}
\label{CONV}
{\rm where} \qquad 
\D(\lambda \mid \{x_k\}) \equiv \Xi(s)\ \ {\rm with\ \ } \lambda \defi s(s-1):
\end{equation}
i.e., $\Xi(s)$ has been rewritten as an infinite product 
$\D(\lambda \mid \{x_k\})$
which is {\sl manifestly even\/} (under $s \longleftrightarrow (1-s)$),
and will qualify as a Delta function of order $\mu_0=\hf$
{\sl in the variable\/} $z=\lambda$, 
a much simpler situation than $\mu_0=1$ (naively suggested by eq.(\ref{CHAD})).
We now derive the ensuing properties of $Z(\sigma)$ 
(as compiled in Sec.7, Table 1).

\subsection{Properties of $Z(\sigma)$ for $\Re \sigma <1$}

As basic initial result, the sequence $\{ x_k \defi \qt +{\tau_k}^2 \}$ 
{\sl is admissible of order\/} $\hf$: 
it obviously fulfills assumptions (i)--(ii);
the function $\Xi(s)$ is entire of order 1 in $s$, hence $\hf$ in $\lambda$;
finally, a large-$\lambda$ expansion of the form (\ref{DAs}) for 
$\log \D (\lambda \mid \{x_k\})$ is easily obtained as follows.
In eq.(\ref{PROD}), for $s \to \infty \ (|\arg s| < \pi/2)$, 
$-\log F(s)$ can be replaced by its complete Stirling expansion  
and $\log \zeta(s) = O(s^{-\infty})$ can be deleted, giving
\begin{equation}
\label{Z1as}
\log \D (\lambda) \sim \log \lambda - {\log \pi \over 2} s 
+ {s \over 2} (\log s -\log 2 -1) - {1 \over 2} \log {s \over 2} 
+{1 \over 2} \log 2 \pi 
+ \sum_{m \ge 1} {B_{2m} 2^{2m-1}\over 2m(2m-1)} s^{-2m+1};
\end{equation}
whereupon $s$ is to be substituted by the relevant solution branch 
of $ s(s-1)=\lambda $, namely
\begin{equation}
\label{SUB}
s = \hf + \sqrt{\lambda + \qt}  \sim
\lambda ^{1/2} + \hf + 
\sum_{n=1}^\infty {2^{-2n}\G(3/2) \over n! \G(-n+3/2)}\lambda ^{1/2-n} .
\end{equation}
The resulting $\lambda$-expansion then has all the required properties, 
with the exponents $\{\mu_n=\hf (1-n)\}$ --- giving the order $\mu_0 = \hf$ --- 
and leading coefficients
\begin{eqnarray}
\label{TAa1}
\tilde a_{1/2} \si=\sf \qt \qquad \qquad \qquad \qquad \qquad
a_{1/2} = - \hf (1 + \log 2\pi) \\
\label{TAa2}
\tilde a_0 \si=\sf \textstyle {7 \over 8}  \qquad \qquad \qquad  \qquad \qquad
a_0 = \qt \log 8\pi \\
\label{TAa3}
\tilde a_{-1/2} \si=\sf \textstyle {1 \over 32}  \qquad \qquad \qquad \qquad \quad
a_{-1/2} = -{1 \over 16} \log 2\pi- {1 \over 48} ,\quad  etc.
\end{eqnarray}
This expansion can be computed to any power $\mu_n$ in principle;
still, a reduced general formula for $a_{\mu_n}$ looks inaccessible this way;
by contrast, all the $ \tilde a_{\mu_n}$ with $\mu_n \ne 0$ arise from 
the single substitution of eq.(\ref{SUB}) into the prefactor 
${s \over 2}$ of $ \log s \ (\sim \hf\log \lambda) $ in eq.(\ref{Z1as}),
giving 
\begin{equation}
\label{TA}
\tilde a_{1/2-m}= {2^{-2m-2}\G(3/2) \over m! \G(-m+3/2)}, 
\qquad \tilde a_{-1-m}=0, \qquad m=0,1,2,\ldots .
\end{equation}

Now, by Sec.2.2$a'$--$a'''$), 
a large-$\lambda$ expansion for $\log \D (\lambda)$ 
translates into explicit properties of $Z(\sigma)$ for $\Re \sigma<1$.
Here, $Z(\sigma)$ gets a double pole at each {\sl half\/}-integer $\hf-m
\ (m \in \mathbb N)$, with principal polar term 
${(-1)^m \over \pi} (\hf -m) \tilde a_{1/2-m}/
(\sigma +m-\hf)^2$, 
and is regular elsewhere; the leading pole $\sigma= \hf$ has 
\begin{equation}
\label{ZPP}
\mbox{full polar part} \ =
{1 \over 8 \pi}{1 \over (\sigma-\hf)^2} - 
{\log 2 \pi \over 4 \pi}{1 \over \sigma-\hf}.
\end{equation}
At $\sigma=0$, eqs.(\ref{TAa2}) and (\ref{ZP0}) 
deliver {\sl two\/} explicit values,
\begin{equation}
\label{Z10}
\textstyle
Z(0) = \tilde a_0 = {7 \over 8}, \qquad \qquad 
Z'(0) = a_0 = \qt \log 8\pi \ (\approx 0.806042857) ;
\end{equation}
the latter makes the {\sl Stirling constant\/}
(the value $\exp-[Z'(0)]$, or regularized determinant) 
also explicitly known for this sequence $\{x_k\}$.

Further quantities, 
tied to the  yet unspecified general coefficients $a_{\mu_n}$,
will also acquire explicit closed forms:
the polar terms of order $1/(\sigma +m-\hf)$
and the finite values $Z(-m)$ for all $m \in \mathbb N$, see Table 1;
but we will need a more indirect approach (Sec.6.1).

\subsection{Properties of $Z(\sigma)$ for $\Re \sigma> \hf $}

Now by Sec.2.2$b'$), $Z(\sigma)$ is holomorphic
in the half-plane $\Re \sigma> \hf $,
and the values $Z(n)$ for $n=1,2,\ldots$ 
lie in the Taylor series of $ \log \D (\lambda) $ at $\lambda=0$,
which can be specified here through eq.(\ref{PROD}).

\medskip

We first expand $ \log \D (\lambda) $ in powers of $s$, where $\lambda=s(s-1)$.
Eq.(\ref{CHAD}) directly implies
\begin{equation}
\label{ZZ}
 \log \D (\lambda)  \equiv \log \Xi(s)=
-\sum_{n=1}^\infty {{\mathscr Z}_n \over n} s^n \qquad 
({\mathscr Z}_n \equiv \sum_\rho \rho^{-n});
\end{equation}
then, from eq.(\ref{PROD}), $\log[\zeta(s)/F(s)]$ and 
$\log[\zeta(1-s)/F(1-s)]$ (expressing the functional equation)
respectively yield the two Taylor series
\begin{eqnarray}
\label{EX1}
-\log \sqrt \pi \ s + \log \G(1+s/2) + \log (1-s) + \log[-2 \zeta(s)]
\qquad \qquad \qquad \qquad \qquad \qquad \nonumber\\
 = -\log \sqrt \pi \, s + \Bigl[ 
-{ \g \over 2} s + \sum_{n=2}^\infty {(-1)^n \zeta(n) \over 2^n n}s^n \Bigr]
- \sum_{n=1}^\infty {1 \over n} s^n 
+ \sum_{n=1}^\infty {(\log |\zeta|)^{(n)}(0) \over n!} s^n , \\ 
\label{EX2}
\log (1-s) + \Bigl[ {\log \pi \over 2} (s-1) +
\log \G \bigl( {1-s \over 2} \bigr) \Bigr] + \log [-s \zeta(1-s)] 
\qquad \qquad \qquad \qquad \nonumber\\
= \Bigl({\g \over 2} + \log 2 \sqrt \pi -1 \Bigr) s +
\sum_{n=2}^\infty {(1-2^{-n}) \zeta(n) - 1 \over n} s^n 
- \sum_{n=1}^\infty { \g^{\rm c}_{n-1} \over (n-1)! } s^n .
\end{eqnarray}
The $\g^{\rm c}_n$ in the last line are {\sl cumulants\/} 
for the Stieltjes constants $ \g_n $ of eq.(\ref{NOT}) \cite{AS} 
(cf. also the $\eta_n$ in ref. \cite{BL}, Sec. 4), i.e., 
\begin{eqnarray}
\label{LAU}
& \displaystyle \log [-s \zeta(1-s)] \equiv
-\sum_{n=1}^\infty {\g^{\rm c}_{n-1} \over (n-1)! } s^n 
\qquad vs \qquad
-s \zeta(1-s) \equiv 
1 - \sum_{n=1}^\infty {\g_{n-1} \over (n-1)! } s^n  \\
&(\g^{\rm c}_0=\g_0=\g, \quad \g^{\rm c}_1 = \g_1 + \hf\g^2 \approx 0.260362078,
\quad \g^{\rm c}_2 = \g_2 + 2\g \g_1 + {2 \over 3} \g^3 \approx 0.034459088,
\ldots).&
\nonumber
\end{eqnarray}
The identification of the three series (\ref{ZZ})--(\ref{EX2})
at each order $s^n$ now yields a countable sequence of 3-term identities:
the first one just restores the result 
${\mathscr Z}_1 = -B$ as in eq.(\ref{Z11});
then, the subsequent ones likewise express the higher ${\mathscr Z}_n$ 
in {\sl two\/} ways,
\begin{equation}
\label{ZRH}
{\mathscr Z}_n =
1 - (-1)^n 2^{-n} \zeta(n) - {(\log |\zeta|)^{(n)}(0) \over (n-1)!}
= 1 - (1-2^{-n}) \zeta(n) + { n \over (n-1)!} \g^{\rm c}_{n-1} ,
\quad n=2,3,\ldots
\end{equation}
That short argument subsumes several earlier results.
The rightmost and center expressions in eq.(\ref{ZRH}) amount to formulae for
${\mathscr Z}_n$ by Matsuoka \cite{Ms2} 
and Lehmer (\cite{L}, eq.(12)) respectively;
the implied relations between 
the derivatives $\zeta^{(n)}(0)$ and the Stieltjes constants $\g_n$
were also discussed in \cite{A,CHO},
together with Euler--Maclaurin formulae for the $\zeta^{(n)}(0)$
which parallel the specification of the $\g_n$ in eq.(\ref{NOT}).

\medskip

As for the values $Z(m)$ themselves, 
they are given by eq.(\ref{ZTay}) now using $\lambda=s(s-1)$ 
as expansion variable, i.e.,
\begin{equation}
\label{ZN}
Z(m) = {(-1)^{m-1} \over (m-1)!} \Bigl[ 
\Bigl( {1 \over 2s-1}{\d \over \d s} \Bigr)^m \log \Xi (s) 
\Bigr] _{s=0 {\rm \ or\ } 1} \qquad (m=1,2,\ldots) .
\end{equation}
But alternatively, $Z(m) \equiv {\mathscr Z}_m + \bigl[ 
\mbox{a {\sl finite\/} linear combination of the } 
\{ {\mathscr Z}_n \}_{n=1,\ldots,m-1} \bigr]$, and vice-versa: as shown below,
\begin{equation}
\label{SZ}
{{\mathscr Z}_n \over n} = \!\! \sum_{0 \le \ell \le n/2} \!\! 
(-1)^\ell {n-\ell \choose \ell} { Z(n-\ell) \over n - \ell} \quad
\Longleftrightarrow \quad 
Z(m) = \sum_{\ell=0}^{m-1} {m+\ell-1 \choose m-1} {\mathscr Z}_{m-\ell}
\end{equation}
($Z(1)={\mathscr Z}_1,\ Z(2)={\mathscr Z}_2+2 {\mathscr Z}_1, \ldots$).
The ${\mathscr Z}_n$ being already known from eqs.(\ref{Z11}) and (\ref{ZRH}), 
$Z(m)$ then reduces to an explicit affine combination (over the rationals) 
of $B=[\log \Xi]'(0)$, $\zeta(n)$, and 
either $(\log |\zeta|)^{(n)}(0)$ or $\g^{\rm c}_{n-1}$ for $1< n \le m$
(see Table 1), as stated earlier by Matiyasevich \cite{MA}. 
(More recently, $Z(1)$ and $Z(2)$ also got revived 
in studies of the distribution of primes \cite{RS,Lb}).

Since ref.\cite{MA} uses eq.(\ref{SZ}) for $Z(m)$ without mentioning any proof,
we sketch one.
First, if $x \defi \rho(1-\rho)$, the expansion of the identity
$\log \bigl[ (1-s/\rho)(1-s/(1-\rho)) \bigr] \equiv \log \bigl[ 1-s(1-s)/x \bigr] $
in powers of $s$ yields 
$ [ \rho^{-n}+(1-\rho)^{-n}] /n  \equiv 
\sum_{0 \le \ell \le n/2} (-1)^\ell {n-\ell \choose \ell} x^{-n+\ell}/(n-\ell)$
for $n=1,2,\ldots$. By recursion, this must invert in the form 
$x^{-m} \equiv \sum\limits_{n=1}^m V_{m,n} [ \rho^{-n}+(1-\rho)^{-n}] $;
then, $\sum\limits_{n=1}^m V_{m,n} \rho^{-n}$ has to be the singular part 
in the Laurent series of $x^{-m}$ around $\rho=0$, 
resulting in $V_{m,n}={2m-n-1 \choose m-1}$. 
Now, summing both sets of identities over the Riemann zeros $\{ \rho \}$ 
yields the stated decompositions (\ref{SZ}).
(We stress that their finite character is specific to the $Z(m)$ as opposed
to all other values $\Z(m,v)$, $v \ne \qt$.)
\medskip

{\sl Note added in proof.\/} ---
For completeness, we quote two other sets of identities 
for the sums ${\mathscr Z}_n$ \cite{VZ}:
$$ 2 {\mathscr Z}_k = - \sum_{\ell=k+1}^\infty 
{\ell-1 \choose k-1} {\mathscr Z}_\ell \qquad \qquad 
\mbox{for each odd } k \ge 1 $$
(a countable sequence of ``sum rules", easy but unreported,
which allow to eliminate any {\sl finite\/} subset of odd values); 
and the connection to {\sl Li's coefficients\/} (cf. [39], thm 2),
$ \lambda_n {\stackrel{\rm def}{=}} \sum_\rho \, [1-(1-1/\rho)^n] $
[which allow to recast 
the Riemann Hypothesis as $ \lambda_n >0 \ (\forall n)$]:
$$ \lambda_n = \sum_{j=1}^n (-1)^{j+1} {n \choose j} {\mathscr Z}_j
\qquad \Longleftrightarrow \qquad 
{\mathscr Z}_n = \sum_{j=1}^n (-1)^{j+1} {n \choose j} \lambda_j \qquad
(n=1,2,\ldots) . $$

\section{Generalized Zeta functions 
$\Z(\sigma,v) \equiv Z(\sigma \mid \{ {\tau_k}^2+v \} )$}

We begin to discuss the Hurwitz-like generalizations of the preceding case
obtained by shifting the {\sl squared\/} parameters ${\tau_k}^2$.
The obviously allowed translations
($x_k \mapsto x_k+v'',\quad  v''>-x_1$ in the real case) preserve the notion of 
admissible sequences (with their values of $\mu_0,\ r$);
that validates the earlier definition
(\ref{ZDef}), as $\Z(\sigma,v) = Z(\sigma \mid \{ {\tau_k}^2+v \} )$.

The corresponding transformation of Delta functions 
(as Hadamard products of order $<1$)
only involves an explicit constant denominator to preserve 
their specific normalization $\D(0)=1$, as
\begin{equation}
\label{SHD}
\D(\lambda \mid \{ {\tau_k}^2+v \} ) \equiv 
\D(\lambda + v-v' \mid \{ {\tau_k}^2+v' \} ) \big/
\D(v-v' \mid \{ {\tau_k}^2+v' \} ).
\end{equation}
Alternatively, we might have opted to normalize Delta functions as
{\sl zeta-regularized\/} infinite products, i.e., 
\begin{equation}
\label{DRG}
\D_{\rm zr}(\lambda) \defi
\exp [-\partial_\sigma Z(\sigma \mid \{ {\tau_k}^2+\lambda \} )]_{\sigma=0} \, ,
\end{equation}
which are fully translation-covariant, but at the same time less explicit.
The overall benefit of this normalization is then dubious 
within the restricted scope of this work; but here,
it explains a dichotomy between algebraic and transcendental properties 
of Zeta functions, which roughly follows our overall division between 
$\{ \Re \sigma <1 \}$ and $\{ \Re \sigma > \mu_0 \}$ properties, but not quite.

\bigskip

Covariance implies that if we translate $\lambda \mapsto (\lambda + v)$,
the expansion of $\log \D_{\rm zr}(\lambda)$ 
around the {\sl invariant\/} point $\lambda=+\infty$ 
can be recomputed to any order by straight substitution,
yielding {\sl explicit polynomials\/} in $v$ as coefficients.
When $\mu_0<1$ as here, then 
$ \D_{\rm zr}(\lambda) \equiv \e^{-a_0} \D(\lambda)$ \cite{V},
hence the previous statement holds for the expansion (\ref{DAs}) 
{\sl minus its term of order\/} $\lambda^0$;
i.e., all the shifted coefficients $\tilde a_{\mu_n}(v), \ a_{\mu_n}(v)$ 
will be polynomial excepting $a_0(v)$.
For $\log \D(\lambda \mid \{ {\tau_k}^2+v \} )$ specifically, 
eqs.(\ref{TAa1}--\ref{TA}) (at $v=\qt$) imply that
\begin{equation}
\label{AA}
\tilde a_{1/2},\ a_{1/2},\ \tilde a_0  \mbox{ stay constant \quad (as well as }
\tilde a_{-1} \equiv \tilde a_{-2} \equiv \cdots \equiv 0);
\end{equation}
\begin{equation}
\label{TAG}
\tilde a_{1/2-m}(v) \equiv {\G(3/2) \over \G(-m+3/2)}{v^m \over 4\,m!}.
\end{equation}
For functions like $\Z(\sigma,v)$, the consequences are that their polar parts 
and ``trace identities" will depend polynomially on $v$;
furthermore, a single (fixed-$v$) large-$\lambda$ expansion, 
such as eq.(\ref{Z1as}) for 
$\log \D(\lambda \mid \{ {\tau_k}^2+{1 \over 4} \} )$, 
suffices to express those $v$-dependences in full. 
Precisely here, by eq.(\ref{AA}): $\Z(\sigma ,v)$ keeps its rightmost 
($\sigma=\hf$) full polar part {\sl constant\/} (and given by eq.(\ref{ZPP})),
as well as its value at 0 ($\Z(0,v) \equiv {7 \over 8}$);
all its other poles (of order 2, except at $v=0$) keep fixed locations.
(As a by-product, any difference function 
$[\Z(\sigma ,v)-\Z(\sigma ,v_0)]$ is holomorphic for $\Re \sigma > -\hf$).
We can specify such polynomial formulae in the half-plane $\Re \sigma < \hf$
still further, but only later by a different path (Sec.6).

\bigskip

By contrast, all formulae for $\Z(\sigma,v)$ 
in the half-plane $\{\Re \sigma > \hf \}$ refer to Taylor coefficients
of $\log \D_{\rm zr}(v+\lambda)$ around  $v$ {\sl finite\/},
which evolve {\sl transcendentally\/} with $v$: here they will only express 
in terms of $\log |\zeta(s)|$ (or $\log \Xi(s)$, from eq.(\ref{XI}))
and its derivatives at $s= \hf \pm  v^{1/2}$.
The first of those coefficients, $\log \D_{\rm zr}(v)$ itself ($= -a_0(v)$),
actually yields a special value lying at $\sigma=0$, by eq.(\ref{DRG}):
\begin{eqnarray}
\label{DETG}
\partial_\sigma  \Z(\sigma,v)_{\sigma=0} = a_0(v) 
\si=\sf a_0(\qt)-\log \D(v- \qt \mid \{ {\tau_k}^2 + \qt \} ) \nonumber\\
\si=\sf \qt \log 8\pi - \log \Xi(\hf \pm v^{1/2})
\end{eqnarray}
(a result which fully matches eq.(\ref{DZet1}) below 
for the Hurwitz-type function $\xi$ \cite{Dn,SS}).
Then, the Taylor coefficients of order $n \ge 1$
(identical for $\log \D(\lambda)$ and $\log \D_{\rm zr}(\lambda)$) yield
\begin{eqnarray}
\label{ZNV}
& \displaystyle \Z(n,v) = {(-1)^{n-1} \over (n-1)!} 
{\d^n \over \d v^n} \log \Xi \bigl( \hf \pm v^{1/2} \bigr) 
\qquad (n=1,2,\ldots), & \\
\label{Z1V}
 e.g., \qquad 
& \Z(1,v) = \pm \hf v^{-1/2} (\log \Xi)'(s=\hf \pm v^{1/2}) \quad (v \ne 0),
\qquad \Z(1,0) = \hf (\log \Xi)''(\hf) & \nonumber
\end{eqnarray}
(but we are in lack of more reduced closed forms for general $(n,v)$).

\medskip

In summary, the polar parts of $\Z(\sigma,v)$ and the special values 
$\{ \Z(-n,v) \}_{n \in \mathbb N}$ (``trace identities") 
have polynomial expressions in $v$; 
whereas $\{ \Z(n,v)\ (n \ne 0$) {\sl plus\/ $\partial_\sigma \Z(0,v)$ at $n=0 \}$}
are also special values, but only computable transcendentally. 
(This conclusion is moreover fully obeyed for typical spectral zeta functions.)

\section{Delta function based at $s= \hf$, and the Zeta function 
$\Z(\sigma) \defi \Z(\sigma,v=0)$}

An interesting option is now to shift the parameter $v$ 
from its initial value $\qt$ in $Z(\sigma)$, 
to the most symmetrical value $v=0$. 
By eq.(\ref{SHD}), the Hadamard product (\ref{CONV}) becomes based at $s=\hf$ as
\begin{equation}
\label{PRD}
\Xi(\hf + t) \equiv \Xi(\hf) \D (t^2 \mid \{{\tau_k}^2\}) \qquad
\mbox{with} \quad t \defi s - \hf , \quad t^2 = \lambda + \qt,
\end{equation}
\begin{equation}
\label{XIV}
\mbox{and} \quad
\Xi(\hf) = - \pi^{-1/4} \G(\fq) \zeta(\hf) \quad (\approx 0.994241556)
\end{equation}
(but very little is known about $\zeta(\hf)$ \cite{Ms1}
and we cannot make this constant factor any more explicit,
contrary to the special case $v=\qt$ where that factor was $\Xi(0)=1$).

The factorized representation (\ref{PROD}) then transforms to 
\begin{equation}
\label{FAC}
\D (t^2 \mid \{{\tau_k}^2\}) {{\bf D}(t) \over 1-2t} = \zeta(\hf + t), \qquad
{\bf D}(t) \defi 
\zeta(\hf) \G(\fq) \pi^{t/2} \big/ \G(\fq + \textstyle{t \over 2}).
\end{equation}
(This Delta function is closest to the determinant of Riemann zeros used
by Berry--Keating for other purposes \cite{BK}.)

We accordingly switch to the Zeta function of the sequence $\{{\tau_k}^2\}$
\cite{G1,Dl}, or in short,
\begin{equation}
\label{Z2}
\Z(\sigma) \defi \Z(\sigma,0)  
= \sum_{k=1}^\infty {\tau_k}^{-2\sigma} , \qquad \Re \sigma > \hf .
\end{equation}
Numerically, this new function looks almost indistinguishable 
from $Z(\sigma)$ (see App.B).
(Also, by Sec.4, $(Z-\Z)(\sigma)$ extends holomorphically to 
$\Re \sigma > - \hf$.)
By contrast, the {\sl meromorphic continuation\/} of $ \Z(\sigma) $ 
will prove to be distinctly simpler than that of $Z(\sigma)$,
thanks to specially explicit representation formulae 
in the half-plane $\{\Re \sigma < \hf \}$.
To obtain these, we now switch to a more powerful, special to $v=0$, approach
(whereas the earlier considerations would still describe $\Z(\sigma)$, 
but just to the same extent as $Z(\sigma)$).

\subsection{The shifted spectrum of trivial zeros}

The factor ${\bf D}(t)$ in eq.(\ref{FAC}) has the structure of a 
{\sl spectral determinant\/} built over the ``spectrum" of trivial zeros
{\sl in the variable\/} $(-t)$, namely $\{\hf+2k\}$
(${\bf D}(t)$ is not exactly the zeta-regularized determinant,
but again this will not matter here).
That spectral interpretation can be extended to the factor $(1-2t)^{-1}$, 
by treating the pole $t=\hf$ (of $\zeta(\hf + t)$)  
as a ``ghost eigenvalue" of multiplicity ($-1$). 
A major role of the spectrum of trivial zeros is to make 
$\log [{\bf D}(t)/(1-2t)]$ asymptotically cancel  
$\log \D(t^2 \mid \{{\tau_k}^2\})$ {\sl to all orders\/} when 
$t \to \infty$ in $|\arg t| < \pi/2$,
given that $\log \zeta(\hf + t)$ decreases exponentially there.

We therefore expect the {\sl spectral zeta function of the trivial zeros\/} 
(of $\zeta(\hf - t)$) to play an important role;
this ``shadow zeta function of $\zeta(s)$" (for short) involves
{\sl both $\zeta(s)$ itself and the partner function} $\beta(s)$,
in the combination
\begin{equation}
\label{SDef}
{\bf Z}(s) \defi \sum_{k=1}^\infty (\hf + 2k)^{-s} \equiv 
2^s \Bigl[\hf \Bigl( (1-2^{-s}) \zeta(s) +\beta(s) \Bigr) -1 \Bigr]
\end{equation}
($= 2^{-s} \zeta(s, \fq)$ in terms of the Hurwitz zeta function).
${\bf Z}(s)$ has a single simple pole at $s=1$, of residue $\hf$, 
and admits the special values 
\begin{eqnarray}
\label{SVS1}
{\bf Z}(-n)
\si=\sf -{2^n \over n+1} B_{n+1}(\qt)-2^{-n} 
= \hf \Bigl[ (1-2^{-n}) {B_{n+1} \over n+1} + 2^{-n-1} E_n \Bigr] - 2^{-n}
, \quad n=0,1,\ldots \\
\label{SVS2}
{\bf Z}(n) \si=\sf {(-1)^n 2^{-n} \over (n-1)!} [\log \G]^{(n)}(\fq), 
\qquad \qquad \qquad \qquad \qquad \qquad  \qquad \qquad \qquad n=2,3,\ldots \\
\si=\sf \left\{  
{ \hf \Bigl[ (2^{2m}-1) {(2\pi)^{2m} \over 2 (2m)!} |B_{2m}| 
+ 2^{2m} \beta(2m) \Bigr] -2^{2m}, \quad n=2m  \atop
 \hf \Bigl[ (2^{2m+1}-1)  \zeta(2m+1) 
+ {\pi^{2m+1} \over  2 (2m)!} |E_{2m}|  \Bigr] -2^{2m+1}, \quad n=2m+1} \right\}
m=1,2,\ldots \nonumber \\
\label{SVS3}
{\bf Z}(0) \si=\sf 
\textstyle -{3 \over 4} \qquad \mbox{and} \qquad {\bf Z}'(0) =  
-{7 \over 4} \log 2 - \hf \log \pi + \log \G \bigl( {1 \over 4} \bigr).
\end{eqnarray}

Remark: literally, the framework of Sec.2 excludes the sequence of trivial zeros
(of linear growth, and order $\mu_0=1$), 
but the truly relevant function here will be ${\bf Z}(2 \sigma)$, 
as Zeta function of the modified sequence $\{ (\hf + 2k)^2 \}$,
which is admissible of order $\hf$ again.

\subsection{Meromorphic continuation formulae for $\Z(\sigma)$}

We start from a slight variant of the representation (\ref{ZMel}) 
for $\Z(\sigma)$, obtained through an integration by parts upon 
the Mellin formula (\ref{IMel}) 
(where $z \equiv t^2$, by eq.(\ref{PRD})):
\begin{equation}
\label{JMel}
\Z(\sigma) = {\sin \pi \sigma \over \pi} J(\sigma), \qquad
J(\sigma) \defi \int_0^{+\infty} (t^2)^{-\sigma} 
\d\,\log \D(t^2 \mid \{{\tau_k}^2\}) \quad 
(\hf < \Re \sigma <1) .
\end{equation}
We next introduce a (regularized) {\sl resolvent trace\/} 
for the spectrum of trivial zeros,
\begin{equation}
\label{RS}
{\bf R}(t) \defi {\d \over \d t} \log {\bf D} (t) = \hf \Bigl[ \log \pi 
- { \G' \over \G} \bigl(\fq + \textstyle{t \over 2} \bigr) \Bigr],
\end{equation}
which has a simple pole of residue $+1$ at each trivial zero of 
$\zeta (\hf+t)$; 
a corresponding function for the pole (``ghost") at $t={1\over 2}$ is 
\begin{equation}
\label{RG}
R_{\rm g}(t) = -1 / (t - \hf) \qquad (\mbox{with residue } (-1)).
\end{equation}
Then, upon insertion of the factorization formula (\ref{FAC}), 
eq.(\ref{JMel}) yields
\begin{equation}
\label{JV}
J(\sigma) \equiv \int_0^{+\infty} t^{-2\sigma} \Bigl[ -{\bf R}(t) - R_{\rm g}(t)
+ {\zeta ' \over \zeta} \bigl(\hf+t \bigr) \Bigr] \d t .
\end{equation}
Now a crucial feature of the case $v=0$ is that {\sl this\/} integral 
is also a Mellin transform with respect to the argument appearing 
in the factorized form of $\zeta(s)$ 
(namely the variable $t$, in eq.(\ref{FAC})).
As a consequence, the contribution to $J$ from ${\bf R}(t)$ 
(and also $R_{\rm g}(t)$) 
can be neatly extracted and evaluated, in closed and interpretable form.
Because the factor in brackets in eq.(\ref{JV}) is $O(t)$ at $t=0$
(due to the functional equation), $J(\sigma)$ can be split as
\begin{eqnarray}
\label{JSR}
J(\sigma)={\bf J}(\sigma)+J_{\rm r}(\sigma), \qquad {\bf J}(\sigma) \si\defi\sf 
\int_0^{+\infty} t^{-2\sigma}[{\bf R}(0)-{\bf R}(t)] \d t \\
J_{\rm r}(\sigma) \si\defi\sf 
\int_0^{+\infty} t^{-2\sigma} \Bigl[ R_{\rm g}(0)-R_{\rm g}(t) +
{\zeta ' \over \zeta} \bigl(\hf + t \bigr) 
- {\zeta ' \over \zeta} \bigl(\hf \bigr) \Bigr] \d t , \nonumber
\end{eqnarray}
({\sl this\/} splitting preserves the convergence strip 
$\{ \hf < \Re \sigma <1\}$ for both resulting integrals).

$J_{\rm r}(\sigma)$ can be split still further, 
once its integration path has been rotated by a small angle: 
either $+\eps$ or $-\eps$, 
in order to bypass the poles of $R_{\rm g}(t)$ and of
${\zeta ' \over \zeta} \bigl(\hf + t \bigr)$ at $t = \hf$,
\begin{eqnarray}
\label{JGZ}
J_{\rm r}(\sigma) = J_{\rm g}^\pm (\sigma) + J_\zeta^\pm (\sigma), \qquad  
J_{\rm g}^\pm (\sigma) \si\defi\sf
\int_0^{+\e^{\pm \mi \eps} \!\infty} t^{-2\sigma} [R_{\rm g}(0)-R_{\rm g}(t)] 
\d t \\
 J_\zeta^\pm (\sigma) \si\defi\sf \int_0^{+\e^{\pm \mi \eps} \!\infty} t^{-2\sigma} 
\Bigl[ {\zeta ' \over \zeta} \bigl(\hf + t \bigr) 
- {\zeta ' \over \zeta} \bigl(\hf \bigr) \Bigr] \d t . \nonumber
\end{eqnarray}

Now, ${\bf J}(\sigma),\ J_{\rm g}^\pm (\sigma)$ can be straightforwardly
transformed into Hankel contour integrals and then 
{\sl computed in closed form\/} (by the residue calculus), giving
\begin{equation}
\label{JSG}
{\bf J}(\sigma) = -{\pi\, {\bf Z}(2\sigma) \over \sin 2 \pi \sigma}, \qquad
J_{\rm g}^\pm (\sigma) = {\pi\, 2^{2\sigma} 
\e ^{\mp 2 \pi \mi \sigma} \over \sin 2 \pi \sigma} ,
\end{equation}
both of which are explicit functions, meromorphic in the whole plane;
chiefly, ${\bf J}$ brings in the shadow zeta function (\ref{SDef}).

Then, again upon back-and-forth integrations by parts, 
$J_\zeta^\pm (\sigma)$ continue to
\begin{equation}
\label{JZ}
J_\zeta^\pm (\sigma) \equiv \int_0^{+\e^{\pm \mi \eps} \!\infty} t^{-2\sigma} 
{\zeta ' \over \zeta} \bigl(\hf + t \bigr) \d t \qquad 
\mbox{analytic for } -\infty < \Re \sigma < \hf ,
\end{equation}
and (cf. eq.(\ref{IP1}))
these integrals admit meromorphic extensions to the whole plane, with
\begin{equation}
\label{JP}
\mbox{simple poles at } \sigma={n \over 2}, \quad \mbox{of residues } 
-{1 \over 2} {(\log |\zeta|)^{(n)} (\hf) \over (n-1)!}, \quad
n=1,2,\ldots
\end{equation}
(the difference 
$J_\zeta^+(\sigma) - J_\zeta^-(\sigma)\equiv 2 \mi \pi 2^{2\sigma}$ 
is entire).

\medskip

All in all, we finally get two complex conjugate Mellin representations 
for $ \Z(\sigma) $:
\begin{equation}
\label{ZC}
\Z(\sigma) = 
{-{\bf Z}(2\sigma)+ 2^{2\sigma} \e^{\mp 2 \pi \mi \sigma} \over 2 \cos \pi \sigma}
+ {\sin \pi \sigma \over \pi} \int_0^{+\e^{\pm \mi \eps} \!\infty} t^{-2\sigma}
{\zeta ' \over \zeta} \bigl(\hf + t \bigr) \d t ;
\end{equation}
and one real principal-value integral representation given by their half-sum,
\begin{equation}
\label{ZR1}
\Z(\sigma) = 
{-{\bf Z}(2\sigma)+ 2^{2\sigma} \cos 2 \pi \sigma \over 2 \cos \pi \sigma}
+ {\sin \pi \sigma \over \pi} \ {\bf -} \!\!\!\!\!\! \int_0^{+\infty} t^{-2\sigma}
{\zeta ' \over \zeta} \bigl(\hf + t \bigr) \d t ,
\end{equation}
(each of the above converges in the full half-plane $\{\Re \sigma < \hf \}$).
 
\noindent Another real form can be obtained with a regular integrand, 
directly from eq.(\ref{JSR}): 
\begin{equation}
\label{ZR2}
\Z(\sigma) = 
- {{\bf Z}(2\sigma) \over 2 \cos \pi \sigma}
+ {\sin \pi \sigma \over \pi} \int_0^{+\infty} t^{-2\sigma}
\Bigl[ {\zeta ' \over \zeta} \bigl(\hf + t \bigr)+ {1 \over t - \hf} \Bigr]
\d t ,
\end{equation}
however this last integral only converges in the strip 
$\{0 < \Re \sigma < \hf \}$.

Remarks: 

- as analytical extension formulae, eqs.(\ref{ZC}--\ref{ZR2})
are precise {\sl counterparts of the functional equation\/} for $\zeta(s)$;
they also stand as more explicit forms of Guinand's functional relation 
for $\Z(\sigma)$ \cite{G1}, as discussed below (eq.(\ref{GFE})).

- a similar formula exists for the function $\xi(s,x)$ of eq.(\ref{ZZet})
(\cite{Dn}, middle of p.149), only requiring $\Re x >1$
(which precisely avoids the problem raised above by the pole of $\zeta(s)$);
in comparison, the present results correspond to the fixed value $x=\hf$, 
since $\Z(\sigma) \equiv 
(2\pi)^{-2\sigma}(2\cos \pi \sigma)^{-1} \xi(2\sigma,\hf)$ by eq.(\ref{Z2X});

\medskip

As we will elaborate next,
analytical properties of $\Z(\sigma)$ in the half-plane $\{\Re \sigma < \hf \}$
are made totally straightforward by the Mellin formulae (\ref{ZC}--\ref{ZR2})
(while $\Z(\sigma)$ is holomorphic in the half-plane $\{\Re \sigma > \hf \}$, 
where its defining series (\ref{Z2}) converges).
Detailed results are also recollected in fully reduced form in Sec.7, Table 1.

\subsection{Properties of $\Z(\sigma)$ for $\Re \sigma < 1$}

A few leading properties of $\Z(\sigma)$ in the half-plane 
$\{\Re \sigma < 1 \}$ emerge more easily as special cases from Sec.4
(although they can be drawn from eq.(\ref{ZC}) as well):

- $\sigma =\hf$ is a {\sl double pole\/}, 
with the {\sl same full polar part\/} as for $Z(\sigma)$, eq.(\ref{ZPP});

- by specializing eqs.(\ref{AA}),(\ref{DETG}),
\begin{equation}
\label{Z20}
\textstyle
\Z(0)=Z(0)={7 \over 8} ; \qquad 
\Z'(0) = Z'(0) - \log \Xi \bigl( \hf \bigr) \ (\approx 0.811817944).
\end{equation}

Otherwise, a Mellin representation like (\ref{ZC}) gives a better global view 
of $\Z(\sigma)$ over the whole half-plane $\{\Re \sigma < \hf\}$.
Indeed, its non-elementary part (the integral (\ref{JZ})) 
becomes {\sl regular\/} there,
hence can be {\sl ignored\/} both for the polar analysis and 
(thanks to the $\sin \pi \sigma$ factor) 
for the ``trace identities" at integer $\sigma$: 
all such information lies then in the {\sl first term\/} alone,
accessible by mere inspection. We thus obtain that:

- $\Z(\sigma)$ only has {\sl simple poles\/} 
at the negative {\sl half\/}-integers $\sigma=\hf-m$, with residues
\begin{equation}
\label{Z2P}
{\mathcal R}_m  = {(-1)^m \over 2 \pi} [{\bf Z}(1-2m) + 2^{1-2m}]
\equiv {(-1)^m \over 8 \pi m} (1-2^{1-2m}) B_{2m} , \qquad  m=1,2,\ldots
\end{equation}
(hence, only the leading pole $\sigma=\hf$ stays double);

- at the negative integers, the {\sl ``trace identities"\/} read as 
\begin{equation}
\label{Z2T}
\Z(-m) = {(-1)^m \over 2} [-{\bf Z}(-2m) + 2^{-2m}]
\equiv (-1)^m 2^{-2m} (1 - {1 \over 8} E_{2m}),  \qquad m=0,1,\ldots
\end{equation}
(both formulae (\ref{Z2P}),(\ref{Z2T}) were fully reduced 
using eq.(\ref{SVS1})).

- a ($\sigma \to -\infty$) {\sl asymptotic formula\/} follows 
for $J_\zeta^\pm(\sigma)$,
from the term-by-term substitution of the Euler product for $\zeta(s)$ 
into the integrand of eq.(\ref{ZC}), giving
\begin{equation}
\label{ZAs}
J_\zeta^\pm(\sigma) \sim 
-\G(1-2\sigma) \sum_{n \ge 2} \Lambda(n) n^{-1/2} (\log n)^{2\sigma-1}, 
\quad \sigma \to -\infty
\end{equation}
where as usual, $\Lambda(n) = \log p$ if $n=p^r$ for some prime $p$, else 0.
(An asymptotic formula for $\Z(\sigma)$ itself then follows from
eq.(\ref{ZC}) and the functional equations for $\zeta(s),\ \beta(s)$.)

Remark: in our notations, Guinand's functional relation for $\Z(\sigma)$
\cite{G1} reads as
\begin{equation}
\label{GFE}
\Z(\sigma) = 
- {{\bf Z}(2\sigma) \over 2 \cos \pi \sigma}
- {Z_{\rm p}(1-2 \sigma) \over 2 \G(2 \sigma) \cos \pi \sigma},
\end{equation}
where 
$Z_{\rm p}(1-2 \sigma) \defi \!\!
 \lim\limits_{T \to +\infty} \biggl\{ \displaystyle
\sum_{2 \le n < \e^T} \!\!\!\! \Lambda(n) n^{-1/2} (\log n)^{2\sigma-1}
- \int_0^T  \e^{x/2} x^{2\sigma-1} \d x \biggr\} $ 
--- subject to the Riemann Hypothesis \cite{CHA} --- clearly specifies 
a (real-valued) resummation of the {\sl divergent\/} series in eq.(\ref{ZAs})
(the asymptotic series for $- J_\zeta^\pm(\sigma) /\G(1-2\sigma) $).
Eq.(\ref{GFE}) was only asserted for $0< \Re \sigma < \qt$,
with no clue as to the analytic structure 
of either $\Z(\sigma)$ or $Z_{\rm p}(1-2 \sigma)$ elsewhere.
The present formulae (\ref{ZC}--\ref{ZR2}) are thus 
resummed versions of eq.(\ref{GFE}), with a definitely more explicit content.

\subsection{Properties of $\Z(\sigma)$ for $\Re \sigma > \hf$}

As stated before, $\Z(\sigma)$ is regular in the half-plane 
$\{\Re \sigma >\hf\}$, where analytical results are identities directly
obtainable by expanding the logarithm of the functional relation (\ref{FAC}) 
in Taylor series around $t=0$. 
Here we will extract those results from the Mellin representation (\ref{ZC}),
invoking the meromorphic properties of its integral term in the whole plane 
as given by eq.(\ref{JP}).

- for half-integer $\sigma=\hf +m$: 
the residues of the two summands in (\ref{ZC}) have to cancel 
given that $\Z(\sigma)$ is analytic in the half-plane; 
this imposes 
\begin{eqnarray}
\label{ZID}
(\log |\zeta|)^{(2m+1)} \bigl( \hf \bigr) \si=\sf
(2m)!\ [ {\bf Z}(2m+1) + 2^{2m+1}],\quad  m=1,2,\ldots \\
(\si=\sf -2^{-2m-1} (\log \G)^{(2m+1)} (\qt) \ ), \nonumber
\end{eqnarray}
which simply amounts to $(\log \Xi)^{(2m+1)}\bigl( \hf \bigr) = 0$
(itself a consequence of the functional equation 
$\Xi(\hf+t) \equiv \Xi(\hf-t)$);
that result further reduces, using eq.(\ref{SVS2}), to the identity
\begin{equation}
\label{ZId}
(\log |\zeta|)^{(2m+1)} \bigl( \hf \bigr) =
\hf (2m)! \, (2^{2m+1}-1)  \zeta(2m+1) + \qt \pi^{2m+1} |E_{2m}| .
\end{equation}
The case $m=0$ is singular, but $(\log \Xi)'(\hf)=0$ directly yields
\begin{equation}
\label{ZI1}
{\zeta' \over \zeta}(\hf) = \hf \Bigl[ \log \pi - (\log \G)' (\qt) \Bigr] =
\hf \log 8\pi + {\pi \over 4} +{\g \over 2} \quad (\approx 2.68609171).
\end{equation}

- for integer $\sigma=m$, the pole of the integral is cancelled 
by the zero of $\sin \pi \sigma$, 
and the following explicit relation results,
\begin{equation}
\label{Z2V}
2 (-1)^{m+1} \Z(m) 
- {1 \over (2m-1)!} (\log |\zeta|)^{(2m)} \bigl( \hf \bigr) =
 {\bf Z}(2m) - 2^{2m}, \quad  m=1,2,\ldots.
\end{equation}
which can also be further reduced with the help of eq.(\ref{SVS2}), see Table 1.

Unfortunately, we hardly know anything else about the values 
$(\log |\zeta|)^{(n)} \bigl( \hf \bigr),\ n =0,1,\ldots$ 
(cf. \cite{Ms1} for $\zeta(\hf)$).
To supplement the relation (\ref{ZId}) with $\zeta(n)$ for $n$ odd,
we can only refer to other formulae for $\zeta(2m+1)$ (compiled in \cite{BBC}),
and to Euler--Maclaurin formulae for $\zeta^{(n)}(s)$
(valid at $s=\hf$) with related numerical data \cite{A,CHO,BBC}.
So, even at $v=0$, the transcendental values $\partial_\sigma \Z(0,v)$
(eq.(\ref{Z20})) and $\Z(m,v)$ currently remain more elusive than at
the (exceptional) point $v=\qt$ (Sec.3.3).
Furthermore, we found no reference at all to those values 
(i.e., $\Z'(0)$, $\Z(m)$) in the literature.

\subsection{Speculations and generalizations}

The results of Secs.5.2--4 for $\Z(\sigma)$ are similar to those yielded by
the ``sectorial" trace formula for the analogous spectral zeta function
$\Z_X(\sigma)$ over a compact hyperbolic surface $X$ \cite{CV,Vz}. 
The present formulae for the Riemann case nevertheless show 
several distinctive features.

\noindent - as announced end of Sec.2.1, the sequence $\{ {\tau_k}^2 \}$ 
and the analogous spectrum of the Laplacian on $X$ 
have mutually singular features:
the former has the parameter values $\mu_0=\hf, \ r=2$ 
($\Z(\sigma)$ has its leading pole {\sl double\/}, at $\sigma=\hf$),
whereas the latter more precisely has $\mu_0=1, \ r=1$ 
($\Z_X(\sigma)$ has all its poles {\sl simple\/}, starting at $\sigma=1$),
hence this spectral analogy for the Riemann zeros holds only partially;

\noindent - in the continuation formulae (\ref{ZC}--\ref{ZR2}), 
$\zeta (s)$ itself {\sl reenters\/} 
as an {\sl additive\/} component of the shadow zeta function ${\bf Z}(s)$. 
This is an altogether different incarnation of $\zeta (s)$ from its initial,
{\sl multiplicative\/} involvement, 
which remains in the integral term and indirectly through the zeros,
in the left-hand side $\Z(\sigma)$. 
It is curious to find two such copies of $\zeta (s)$ to coexist in one formula, 
especially with the additive $\zeta (2 \sigma)$ represented 
in its critical strip; 


\noindent - however, those formulae relative to $\zeta(s)$ 
are not fully closed as they also invoke the other Dirichlet series
$\beta (2 \sigma)$ 
(as second additive component in the shadow zeta function ${\bf Z} (2 \sigma)$).
The question then arises whether $\beta(s)$ and other zeta functions
can be handled on the same footing as $\zeta(s)$ (as in \cite{K}),
so we now outline a possible extension of eqs.(\ref{ZC}--\ref{ZR2}). 

\medskip

We assume that $\tilde \zeta (s)$ is a Dirichlet zeta or $L$-series, having:

\noindent - a single pole, at $s=1$ and of order $q$ (typically, $q=0$ or 1); 

\noindent - the asymptotic property $\log \tilde \zeta(s) = o(s^{-N})$ 
for all $N\ (s \to +\infty)$;

\noindent - a functional equation of the form
\begin{equation}
\label{FE}
\tilde {\D} (t^2) {\tilde {\bf D} (t) \over (1-2t)^q} = 
\tilde \zeta(\hf + t), \quad \mbox{where}
\end{equation}
 
\noindent - $\tilde {\D} (t^2)$ is an entire function of order $<1$ 
{\sl in the variable\/} $t^2$, and

\noindent - $\tilde {\bf D} (t)$ is an entire function 
with all of its zeros lying in the half-plane $\{ \Re t < 0 \}$.

Then the Zeta functions $\Z_{\tilde \zeta}(\sigma)$ 
(for the zeros of $\tilde {\D}$) and
${\bf Z}_{\tilde \zeta}(\sigma)$ (for the zeros of $\tilde {\bf D}$) 
are related by this formula corresponding to eq.(\ref{ZR2})
(we omit the others),
\begin{equation}
\label{ZR}
\Z_{\tilde \zeta}(\sigma) = 
-{{\bf Z}_{\tilde \zeta} (2\sigma) \over 2 \cos \pi \sigma}
+ {\sin \pi \sigma \over \pi} \int_0^{+\infty} t^{-2\sigma}
\Bigl[ {\tilde \zeta ' \over \tilde \zeta} \bigl(\hf + t \bigr)
+ {q \over t - \hf} \Bigr] \d t .
\end{equation}

Apart from $ \zeta(s) $ itself, with eq.(\ref{ZR2}), 
the next independent example is $\beta(s)$. 
Its functional equation (\ref{xI}) has the form (\ref{FE}) with $q=0$
(no pole) and $ \tilde {\D} (t^2) = \Xi_{\chi_4} (\hf + t) $,  
$\tilde {\bf D} (t) = \bigl( {\pi \over 4} \bigr) ^{{3 \over 4} + {t \over 2}} 
\big / \G \bigl( {3 \over 4}+ {t \over 2} \bigr) $;
its spectrum of trivial zeros (for $\beta(\hf -t)$) is 
$\{-\hf + 2k \}\quad ( = {3 \over 2},\ {7 \over 2}, \cdots)$,
giving as shadow zeta function 
\begin{equation}
\label{SLD}
{\bf Z}_\beta(\sigma)= \sum_{k=1}^\infty (-\hf +2 k)^{-s} \equiv 
2^s \Bigl[\hf  \bigl((1-2^{-s}) \zeta(s) -\beta(s)\bigr)  \Bigr] .
\end{equation}
Under $q=0$, all Mellin representations (\ref{ZC}--\ref{ZR2}) 
coalesce into the single regular form
\begin{equation}
\label{LR}
\Z_\beta(\sigma) = 
-{{\bf Z}_\beta (2\sigma) \over 2 \cos \pi \sigma}
+ {\sin \pi \sigma \over \pi} \int_0^{+\infty} t^{-2\sigma}
{\beta ' \over \beta} \bigl(\hf + t \bigr) \d t  \qquad (\Re \sigma < \hf),
\end{equation}
and all consequences previously drawn for $\Z(\sigma)$ 
have analogs for $\Z_\beta(\sigma)$.

Various such integral representations will naturally {\sl add\/},
whenever the initial zeta functions combine nicely under multiplication.
For instance, eqs.(\ref{ZR2}) and (\ref{LR}) add up to:
\begin{equation}
\label{ZLR}
(\Z_\zeta + \Z_\beta)(\sigma) = 
-{ 2^{2\sigma} [(1-2^{-2\sigma}) \zeta (2\sigma) -1] \over 2 \cos \pi \sigma}
+ {\sin \pi \sigma \over \pi} \int_0^{+\infty} t^{-2\sigma}
\Bigl[ \Bigl( {\zeta ' \over  \zeta} + {\beta ' \over  \beta} \Bigr)
\bigl(\hf + t \bigr) + {1 \over t - \hf} \Bigr] \d t .
\end{equation}
Here, the shadow zeta function purely invokes $\zeta(s)$;
on the other hand, under the integral sign we now find 
$\bigl[ {(\beta \zeta)'  \over \beta \zeta} \bigr] \bigl( \hf +t \bigr) $
so that the new multiplicative zeta function is $\beta(s) \zeta(s)$, 
also recognized as ${1 \over 4}$ times $Z_4 (s)$, 
the {\sl zeta function of the ring of Gaussian integers\/} $\mathbb Z [\mi]$
\cite{C}. Hence eq.(\ref{ZLR}) becomes
\begin{equation}
\label{Z4R}
\Z_{Z_4}(\sigma) = 
-{ 2^{2\sigma} [(1-2^{-2\sigma}) \zeta (2\sigma) -1] \over 2 \cos \pi \sigma}
+ {\sin \pi \sigma \over \pi} \int_0^{+\infty} t^{-2\sigma}
\Bigl[ {{Z_4}' \over Z_4} \bigl(\hf + t \bigr)
+ {1 \over t - \hf} \Bigr] \d t ;
\end{equation}
thus, to isolate $\zeta(s)$, here in the additive position, 
we again had to allow a different zeta function elsewhere, 
this time $Z_4(s)$ in the multiplicative position.

Likewise, by subtracting eq.(\ref{LR}) from (\ref{ZR2}) instead,
we could get the shadow zeta function to be $\beta(s)$; then the counterpart
of eqs.(\ref{ZId}),(\ref{ZI1}) is a {\sl fully explicit\/} identity, 
\begin{equation}
\label{ZB}
(\log |\zeta|)^{(2m+1)} \bigl( \hf \bigr) -
 (\log \beta)^{(2m+1)} \bigl( \hf \bigr) 
= \hf \pi^{2m+1} |E_{2m}| + \delta_{m,0} \log 2, \quad m=0,1,\ldots ,
\end{equation}
whereas each of the two left-hand-side terms separately needs $\zeta(2m+1)$
(or $\g$ for $m=0$).

\section{More about the Hurwitz-type functions}

The purpose of this Section is twofold. 
First, we analyze the Zeta functions 
$\Z(\sigma,v) = \sum_k ({\tau_k}^2+v)^{-\sigma}$
more explicitly over the half-plane $\{ \Re \sigma <1 \}$ than in Sec.4,
by exploiting the latest special properties of the function $\Z(\sigma,0)$
(with new results even for the case $v=\qt$).
Then, by the same approach, we (briefly) discuss the other 
Hurwitz-type Zeta functions ${\mathfrak Z}(\sigma,a)$ and $\xi(s,x)$,
defined through eqs.(\ref{Hz}) and (\ref{ZZet}) respectively.

\subsection{Further properties of $\Z(\sigma,v)$ for $\Re \sigma <1$}

To describe the Hurwitz-type function $\Z(\sigma,v)$ better,
we now systematically expand it in terms of $\Z(\sigma)$, as
\begin{equation}
\label{SHZ0}
\Z(\sigma,v) = \sum_{k=0}^\infty 
({\tau_k}^2)^{-\sigma} \Bigl( 1 + {v \over {\tau_k}^2} \Bigr)^{-\sigma}
= \sum_{\ell=0}^\infty {\G(1-\sigma) \over \ell !\G(1-\sigma-\ell)}
\Z(\sigma+\ell) \, v^\ell \qquad (|v| < {\tau_1}^2).
\end{equation}
Such an expansion can be formulated around any reference point $v_0$, 
but it will be specially useful for $v_0=0$ as above.
For instance, coupled with eq.(\ref{ZC}) (say), it can express the meromorphic
continuation of the general $\Z(\sigma,v)$ to $\{ \Re \sigma < \hf \}$, 
while we lack an analog of eq.(\ref{ZC}) itself for any $v \ne 0$.

For the polar structure of $\Z(\sigma,v)$ at $\sigma=-m+\hf,\ m \in\mathbb N$, 
the series (\ref{SHZ0}) reduces to
\begin{equation}
\label{HP}
\Z(-m+\hf+\eps,v) = 
\sum_{\ell=0}^m {\G(\hf+m-\eps) \over \ell ! \G(\hf+m-\ell-\eps)}
\Z(-m+\ell+\hf+\eps) \, v^\ell
\quad [\mbox{+regular part for } \eps \to 0];
\end{equation}
then, importing the polar structure of $\Z(\sigma)$ from eq.(\ref{Z2P}), we get
\begin{equation}
\label{Hp}
\Z(-m+\hf+\eps,v) = {1 \over 8 \pi} {\G(m+\hf) \over m! \G(\hf)} v^m \,\eps^{-2}
 + {\mathcal R}_m(v) \,\eps^{-1} + O(1) \qquad(\eps \to 0),
\end{equation}
just by brute-force polar expansion of the right-hand side of eq.(\ref{HP}). 
Here, the polar part of order 2 at every $(-m+\hf)$ is clearly induced 
by the only such part of $\Z(\sigma)$ (at $\sigma=\hf$),
through the term with $\ell =m$ in eq.(\ref{HP}); whereas the residue 
${\mathcal R}_m(v)$ is built from all residues of $\Z(\sigma)$ 
at poles with $\sigma \ge -m+\hf$, as
\begin{equation}
\label{HR}
{\mathcal R}_m(v) = -{\G(m+\hf) \over m! \G(\hf)} 
\Bigl[ {1 \over 4 \pi} \sum_{j=1}^m {1 \over 2j-1}
+ {\log 2\pi \over 4 \pi} \Bigr] v^m 
+ \sum_{j=1}^m{\G(\hf+m) \over (m-j)!\G(\hf+j)}
{\mathcal R}_j \, v^{m-j} ,
\end{equation}
(the residues ${\mathcal R}_j$ of $\Z(\sigma)$ at $-j+\hf$
are known from eq.(\ref{Z2P})).

\noindent Remark: for $m=0$, the full polar part (\ref{ZPP}) at $\sigma=\hf$,
independent of $v$, is recovered.

When $\sigma \in -{\mathbb N}$, the series (\ref{SHZ0}) also terminates, as 
\begin{equation}
\label{TIG}
\Z(-m,v) \equiv 
\sum_{\ell=0}^m {m \choose \ell} \Z(-m+\ell)\,v^\ell \qquad (m \in \mathbb N) ,
\end{equation}
so that explicit ``trace identities" for general $v$ derive from 
those for $v=0$ (eq.(\ref{Z2T})). (For $v=\qt$, this result simplifies further,
see Table 1.)

Remark: in view of eqs.(\ref{ZP2}) and (\ref{ZT}), 
the latter two results now imply a general-$n$ formula 
for the coefficients $a_{(1-n)/2}(v)$ 
in the large-$\lambda$ expansion (\ref{DAs}) 
of $\log \D(\lambda \mid \{ {\tau_k}^2+v \})$.
(Hitherto we had such a formula just at $v=0$, not even at $v=\qt$, 
and knew only the {\sl other\/} coefficients $\tilde a_{(1-n)/2}(v)$ 
for any $v$ and $n$, by eqs.(\ref{AA}--\ref{TAG}).)

Our initial emphasis on the special case $v=\qt$ might now seem misplaced:
why didn't we operate at once from $v=0$~?
In the first place, we saw the case $v=\qt$ arise more readily 
from the standard product representation of $\zeta(s)$. 
But mainly, the case $v=\qt$ also enjoys certain special properties, 
this time with the values $Z'(0)$ and $Z(n)$ (Sec.3.3),
and since these evolve from {\sl transcendental\/} functions of $v$ (Sec.4),
their expansions (\ref{SHZ0}) around $v=0$ are now {\sl infinite\/}.
So, each case $v=0$ and $v=\qt$ has its own exceptional features,
the former in the half-plane $\{ \Re \sigma<1 \}$,
and the latter for $\sigma \in \mathbb N$.

\subsection{The Hurwitz-type functions ${\mathfrak Z}(\sigma,a)$ and $\xi (s,x)$}

The function ${\mathfrak Z}(\sigma,a)$ as defined by eq.(\ref{Hz}) 
is {\sl a priori\/} more singular than $\Z(\sigma,a)$ 
(the sequence $\{ \tau_k \}$ itself has $r=2$ and $\mu_0=1$,
which would require a formalism more elaborate than in Sec.2).
Fortunately, ${\mathfrak Z}(\sigma,a)$ can also be analyzed directly
through its expansion around 
${\mathfrak Z}(\sigma,0) \equiv \Z(\sigma)$, by analogy with eq.(\ref{SHZ0})
(see also \cite{HKW}):
\begin{equation}
\label{SHZ1}
{\mathfrak Z}(\sigma,a) = \sum_{k=0}^\infty 
{\tau_k}^{-2\sigma} \Bigl( 1 + {a \over \tau_k} \Bigr)^{-2\sigma}
= \sum_{\ell=0}^\infty {\G(1-2\sigma) \over \ell ! \G(1-2\sigma-\ell)}
\Z(\sigma + \hf \ell) \, a^\ell  \quad (|a|<\tau_1).
\end{equation}
This formula generates a pole for ${\mathfrak Z}(\sigma,a)$
now at every {\sl half-integer\/} $\hf(1-n),\ n \in\mathbb N$, according to:
\begin{equation}
\label{HP1}
{\mathfrak Z}(\hf(1-n)+\eps,a) = 
\sum_{\ell=0}^n {\G(n-2\eps) \over \ell ! \G(n-\ell-2\eps)}
\Z(\hf(1-n+\ell)+\eps) \, a^\ell \quad + O(\eps) .
\end{equation}
Differences with eq.(\ref{HP}) arise due to the factor
$\G(n-2\eps)/\G(n-\ell-2\eps)$ vanishing whenever $\ell \ge n >0$.
Only the polar part at $\sigma=\hf$ remains the same as for $\Z(\sigma,v)$
(of order $r=2$ and independent of $a$, given by eq.(\ref{ZPP})); 
all other poles $\hf(1-n)$ of ${\mathfrak Z}(\sigma,a)$ are now {\sl simple\/}, 
of residues
\begin{equation}
\label{RES1}
r_n(a) = -{1 \over 4 \pi n} \, a^n +
\sum_{0< 2m \le n} { n-1 \choose 2m-1 } {\mathcal R}_m \, a^{n-2m},
\quad n=1,2,\ldots 
\end{equation}
(again, ${\mathcal R}_m$ is the residue given by eq.(\ref{Z2P})).
In addition, at $\sigma=0$ the $\eps$-expansion of eq.(\ref{HP1}) 
captures the {\sl finite part\/} too:
\begin{equation}
\label{RES0}
r_1(a)=\Res_{\sigma=0} {\mathfrak Z}(\sigma,a) = -{a \over 4 \pi}; 
\quad \mbox{finite part: }
{\rm FP}_{\sigma=0} \,{\mathfrak Z}(\sigma,a) = 
{7 \over 8} + { \log 2\pi \over 2\pi } \, a .
\end{equation}

\medskip

As for the function $\xi(s,x)$,
if we express it by eq.(\ref{Z2X}) in terms of ${\mathfrak Z}(\sigma,a)$,
then we find this combination to be {\sl less singular\/} overall:
by mere substitution of eq.(\ref{HP1}), 
$\xi (s,x)$ shows a {\sl simple\/} pole at $s=1$, of residue $-\pi$ \cite{SS},
and all other possible poles at $s=1-n,\ n=1,2,\ldots $ {\sl cancel out\/}, 
resulting in the holomorphy of $\xi (s,x)$ for all $s \ne 1$ 
with the computable finite values (``trace identities")
\begin{equation}
\label{TI1}
\xi (1-n,\hf+y) = {2 \over (2\pi)^{n-1}} \Biggl[ -\pi { r_n( \mi y) \over \mi^n}
+ \sum_{0 \le 2m < n} (-1)^m { n-1 \choose 2m } \Z(-m) y^{n-2m-1} \Biggr] ,
\quad n=1,2,\ldots 
\end{equation}
(An alternative evaluation follows from Deninger's continuation formula 
for $\xi(s,x)$ (\cite{Dn}, middle of p.149), as
$\xi (1-n,\hf+y) = (2\pi)^{1-n} 
\bigl[ (y+\hf)^{n-1} + (y-\hf)^{n-1} + 2^{n-1} B_n(\qt+ \hf y) \bigr] $,
whose agreement with eq.(\ref{TI1}) can be verified.)

As for special values: first, $\partial_s \xi (s,x)_{s=0}$ 
is expressible as well, in terms of $\zeta(x)$ \cite{Dn,SS}:
\begin{eqnarray}
\label{DZet1}
\partial_s \xi (s,x)_{s=0} \si=\sf \log \,2^{1/2} (2\pi)^2 - \log \Xi(x) \\
\label{DZet2}
\Longleftrightarrow \qquad 
-\partial_s \bigl[ \sum_\rho (x-\rho)^{-s} \bigr]_{s=0}
\si=\sf \log \Xi(x) + \hf (\log 2 \pi) \, x - \hf \log 4 \pi
\end{eqnarray}
(the equivalence of the two forms follows from eqs.(\ref{ZZet}) 
and (\ref{TI1}) for $n=1$, i.e., $\xi (0,x) = \hf (x+3)$).
Now, the exponentiated left-hand side of eq.(\ref{DZet2}) precisely defines  
the zeta-regularized product $\tilde {\D}_{\rm zr}(x)$
built upon the sequence $\{\rho\}$ of Riemann zeros,
while the right-hand side mainly involves $\Xi(x)$ of eq.(\ref{CHAD}).
So, eq.(\ref{DZet2}) is converting a zeta-regularized product
($\tilde {\D}_{\rm zr}(x)$) to Hadamard product form.
As an aside, we now verify that such a conversion formula is entirely 
{\sl fixed by universal rules\/} for (complex) {\sl admissible sequences\/},
specialized here to $r=1$ (as in \cite{V,QHS}) and $\mu_0=1$ ---
since the Zeta functions $\xi(s,x)$ have just a simple pole at $s=1$.
Those rules yield these two prescriptions: 
$\log \tilde {\D}_{\rm zr}(x) \equiv \log \Xi(x)-(\alpha x+\beta)$,
{\sl and\/} the large-$x$ expansion of $\log \tilde {\D}_{\rm zr}(x)$
shall only retain {\sl canonical\/} (or {\sl standard\/}) terms, namely: 
$c_\mu x^\mu$ for $1>\mu \notin \mathbb N,\ c_1 x(\log x -1),\ c_0 \log x$.
Those conditions together fix $(\alpha,\beta)$ uniquely, and here,
eqs.(\ref{CONV}--\ref{Z1as}) for $\log \Xi(x)$ as input 
precisely lead to eq.(\ref{DZet2}) as output.

Likewise, the special values $\xi (n,x),\ n=1,2,\ldots $ are expressible
in terms of [the higher Laurent coefficients of] $\zeta(x)$, e.g.,
by applying residue calculus to Deninger's continuation formula 
(\cite{Dn}, p.149). 

Thus, a fair degree of structural parallelism finally shows up
between the two Hurwitz-like families $\Z(s,v)$ and $\xi (s,x)$ \cite{VZ}.

\section{Recapitulation of main results}

In way of conclusion, Table 1 collates the analytical results found
for the two $\zeta$-Zeta functions $ Z(\sigma) \ (=\Z(\sigma,v=\qt))$ 
and $ \Z(\sigma) \ (=\Z(\sigma,0))$. Furthermore, corresponding results 
for the general $\Z(\sigma,v)$ were derived in Sec.4 (for $\Re \sigma \ge 0$) 
and 6.1 (for $\Re \sigma < 1$), 
and partly extended to the functions ${\mathfrak Z}(\sigma,a)$ and $\xi (s,x)$
in Sec.6.2.

\begin{table}
\begin{tabular}  {ccc}
\hline \\[-12pt]
$\sigma$ & 
$Z(\sigma) = \sum\limits_{k=1}^\infty ({\tau_k}^2+\qt)^{-\sigma} \quad [v=\qt]$
& $\Z(\sigma) = \sum\limits_{k=1}^\infty {\tau_k}^{-2\sigma} \quad [v=0]$
\\[10pt]
\hline \\[-12pt]
$-m$ & $ (-1)^{m+1} 2^{-2m-3} 
\sum\limits_{\ell=0}^m {m \choose \ell} (-1)^{\ell}  {E_{2(m-\ell)}}
^{\ {\scriptsize(\ref{TIG},\ref{Z2T})}} $ & 
$ (-1)^m 2^{-2m} (1 - {1 \over 8} E_{2m})\ ^{{\scriptsize(\ref{Z2T})}} $ \\[10pt]
$-m+\hf +\eps$ & 
$ \!\! [ {2^{-2m} \over 8 \pi} {\G(m+1/2) \over m! \G(1/2)}] \,\eps^{-2}
 + {\mathcal R}_m(\qt) \, \eps^{-1} +O (1) \, ^{{\scriptsize(\ref{Hp},\ref{HR})}} $ & 
$ [{(-1)^m \over 8 \pi m} (1-2^{1-2m}) B_{2m}] \,\eps^{-1} + O(1) 
\ ^{{\scriptsize(\ref{Z2P})}} $ \\[3pt]
$\vdots$ & $\vdots$ & $\vdots$ \\[3pt]
$-1$ & $-1/16$ & $-9/32$ \\[6pt]
$-\hf +\eps$ & 
${1 \over 64 \pi} \,\eps^{-2} - [{3 \log 2 \pi + 4 \over 96 \pi}] \,\eps^{-1} 
+O(1)$ &
$-{1 \over 96 \pi} \,\eps^{-1} +O(1)$ \\[6pt]
$0$ & $7/8 \ ^{{\scriptsize(\ref{Z10})}} $ & $7/8 \ ^{{\scriptsize(\ref{Z20})}} $ \\[3pt]
\hline \\[-12pt]
${\mbox{\sl derivative} \atop \mbox{\sl at 0}}$ & 
$Z'(0) = \qt \log 8\pi \ ^{{\scriptsize(\ref{Z10})}} $ &  
$\Z'(0) = \log \bigl[ 2^{11/4} \pi^{1/2} \G (\qt)^{-1} |\zeta (\hf)|^{-1} \bigr]
\ ^{{\scriptsize(\ref{Z20})}} $ \\[6pt]
$+\hf +\eps$ & 
${1 \over 8 \pi} \,\eps^{-2} - {\log 2 \pi \over 4 \pi} \,\eps^{-1} + O(1)
\ ^{{\scriptsize(\ref{ZPP})}} $ & 
${1 \over 8 \pi} \,\eps^{-2} - {\log 2 \pi \over 4 \pi} \,\eps^{-1} + O(1)
\ ^{{\scriptsize(\ref{ZPP},\ref{AA})}} $ \\[6pt]
$+1$ & $-\hf \log 4\pi + 1 + \hf \g \equiv {\mathscr Z}_1 \ ^{{\scriptsize(\ref{Z11})}} $ & 
$ \hf (\log |\zeta|)''(\hf) + {1 \over 8}\pi^2 +\beta(2) -4$ \\[4pt]
$+2$ & 
$\!\!\!\!\!\!\!\!\!\!\!\!\!\!\!\!\!\!
\biggl\{ \! {\textstyle -\log 4\pi+3+\g -(\log |\zeta|)''(0)-{1 \over 24} \pi^2
\atop \textstyle \!\!\!\!\!\!\!\!\!\!\!
-\log 4\pi+3+\g + \, 2\g_1 + \g^2 -{1 \over 8}\pi^2 } \biggr\} $ &
$ -{1 \over 12} (\log |\zeta|)^{(4)} (\hf)  - {1 \over 24}\pi^4 
- 4 \beta(4) + 16 $ \\
$\vdots$ & $\vdots$ & $\vdots$ \\
$m$ & $ \sum\limits_{\ell=0}^{m-1} {m+\ell-1 \choose m-1} {{\mathscr Z}_{m-\ell}} 
^{\ {\scriptsize(\ref{SZ})}}$ &
$\!\!\!\!\!\!\!\!\!\!\!\!\!\!\!\!\!\!\!\!\!\!\!\!\!\!\!\!\!\!\!\!\!\!\!\!\!\!
\matrix{ (-1)^m \Bigl\{ -{1 \over 2(2m-1)!} (\log |\zeta|)^{(2m)} (\hf) \cr
\qquad\qquad\qquad \quad
{} - \qt \bigl[ (2^{2m}-1) \zeta (2m) + 2^{2m} \beta (2m) \bigr] \cr
\ \qquad \qquad \qquad \qquad \qquad \qquad \qquad \qquad \qquad \qquad 
{} + 2^{2m} \Bigr\} ^{{\scriptsize(\ref{Z2V})}} } $ \\[-13pt]
 & \multicolumn{2}{l}{ $\!\!\!\!\!\!\!\!\!\!\!\!\!\!\!\!\!\!\!
\Bigl[ {\mathscr Z}_n \equiv \Biggl\{ {\!\textstyle 
1 - (-1)^n 2^{-n} \zeta(n) - {(\log |\zeta|)^{(n)}(0) \over (n-1)!}
\atop \textstyle \!
1 - (1-2^{-n}) \zeta(n) + { n \over (n-1)!} \g^{\rm c}_{n-1} } \Biggr\}
\; (n=2,3,\ldots) \Bigr] ^{\, {\scriptsize(\ref{ZRH})}} $ } \\[12pt]
\hline
\end{tabular}
\caption{Analytical results for $\zeta$-Zeta functions of the Riemann zeros.
Notations: see eqs.(\ref{ZER}),(\ref{NOT}), (\ref{LDef}) [for $\beta(s)$],
(\ref{LAU}) [for $\g^{\rm c}_n$]; ${\mathscr Z}_m \equiv \sum_\rho \rho^{-m}$;
$m$ stands for any {\sl positive\/} integer, and $\eps \to 0$. 
As an indexing tool,
the superscripts refer to the relevant equation numbers in the main text.
[In the very last formula (for $\Z(m)$), one can still use 
eq.(\ref{SVR}) for $\zeta(2m)$]. 
}
\end{table}

The other novel results we have developed here concern $\Z(\sigma)$
in the half-plane $\{ \Re \sigma < \hf \}$:
the analytical continuation formulae (\ref{ZC}--\ref{ZR2}),
and the $\sigma \to -\infty$ asymptotic formula (\ref{ZAs}) as corollary.
We also came across two elementary (but unfamiliar to us) formulae 
concerning $\zeta(s)$ itself: eqs.(\ref{ZId}),(\ref{ZB}).

Appendix B gives information on some numerical aspects 
based on our use of the 100,000 first Riemann zeros 
(made freely available on the Web by A.M. Odlyzko \cite{O}, 
to whom we express our gratitude). 
We also wish to thank C. Deninger, J.P. Keating, P. Leb\oe uf, V. Maillot,
C. Soul\'e, and the Referee, 
for helpful references and comments.

\section*{Appendix A: Meromorphic Mellin transforms}

We briefly recall the meromorphic continuation argument for a Mellin transform
like eq.(\ref{IMel}), 
$I(\sigma) \defi \int_0^{+\infty} L(z) z^{-\sigma-1} \d z$,
assuming the function $L(z)$ to be regular on $\mathbb R^+$ (for simplicity), with 
\begin{equation}
\label{LAs}
L(z)=O(z^{\nu_0})  \quad (z \to 0^+), \mbox{ and} \quad
L(z) \sim \sum_{n=0}^\infty 
(\tilde a_{\mu_n} \log z + a_{\mu_n})  z^{\mu_n} \quad (z \to +\infty)
\end{equation}
as in eq.(\ref{DAs}) (asymptotic estimates are repeatedly differentiable); 
and crucially, $\mu_0 < \nu_0$.

\medskip

Sequential directed integrations by parts can be used 
(see \cite{J,V,C} for details).

Step 1: $I(\sigma)$ converges for $\mu_0 < \Re z < \nu_0$, and in that strip,
\begin{equation}
\label{LI1}
I(\sigma) =  \int_0^{+\infty} [L(z) z^{-\mu_0}]' 
{ z^{\mu_0 -\sigma} \over \sigma-\mu_0 } \, \d z .
\end{equation}
If $\tilde a_{\mu_0}=0$, this suffices:
the new integral actually converges for  $\mu_1 < \Re z < \nu_0$
(thanks to $z^{\mu_0 +1}[L(z) z^{-\mu_0}]' =O( z^{\mu_1} \log z)$ 
for $z \to \infty$),
hence $I(\sigma)$ is manifestly meromorphic in that wider strip, 
with a simple pole at $\sigma=\mu_0$ of residue 
\begin{equation}
\int_0^{+\infty} [L(z) z^{-\mu_0}]' \, \d z = a_{\mu_0} ;
\end{equation}
furthermore, in the complementary strip $\{ \mu_1 < \Re z < \mu_0 \}$,
backward integration by parts now yields
\begin{equation}
I(\sigma) = \int_0^{+\infty} [L(z) - a_{\mu_0}z^{\mu_0})] z^{-\sigma-1} \d z .
\end{equation}
Then the whole argument can be restarted from here, 
to extend $I(\sigma)$ further (across $\{ \Re \sigma = \mu_1 \}$), and so on:
the case $r=1$ thus gets settled.

Step 2: if $\tilde a_{\mu_0} \ne 0$, 
one more integration by parts upon eq.(\ref{LI1}) yields
\begin{equation}
I(\sigma) = \int_0^{+\infty} [z [L(z) z^{-\mu_0}]' ]'
{z^{\mu_0 -\sigma} \over (\sigma-\mu_0)^2}\,\d z \qquad (\mu_1 < \Re z < \nu_0).
\end{equation}
All previous arguments then carry over, but the pole is now {\sl double\/}, with
\begin{equation}
\mbox{principal polar cefficient:} \quad 
\int_0^{+\infty} [z [L(z) z^{-\mu_0}]' ]' \, \d z = \tilde a_{\mu_0} ,
\end{equation}
and residue $-\int_0^{+\infty} [z [L(z) z^{-\mu_0}]' ]' \log z \, \d z$
(from the residue calculus);
integrations by parts (backwards, and split) reduce the latter to
\begin{equation}
\mbox{residue:} \quad \int_0^1  [L(z) z^{-\mu_0}]' \, \d z + 
\int_1^{+\infty}  \{ z[L(z) z^{-\mu_0}]' - \tilde a_{\mu_0} \}{1 \over z}\, \d z
= a_{\mu_0} .
\end{equation}
The last two formulae thus generate eq.(\ref{IP2}) for the leading double pole,
and so on for $r=2$.

(More generally, if the factor of $z^{\mu_n}$ in the expansion (\ref{LAs})
is a polynomial of degree $p_n$ in $\log z$, 
then $\mu_n$ becomes a pole of order $(p_n+1)$ for $I(\sigma)$ \cite{J}.)

\medskip

As for meromorphic continuation in the other direction 
(across $\{ \Re \sigma = \nu_0 \}$),
it works likewise if $L(z)$ admits a $z \to 0$ expansion:
the previous arguments apply upon exchanging the bounds 0 and $+\infty$
under $\sigma \mapsto -\sigma$.
E.g., in the regular case, $L(z)$ expands in an entire series at $z=0$, 
and step 1 suffices (as in the main text, where $\nu_0=1$).

\section*{Appendix B: Numerical aspects}

We complete our analytical study by describing some {\sl very heuristic\/}
numerical work with $Z(\sigma )$ and $\Z(\sigma )$,
mostly in the range $\{ \sigma \ge 0 \}$, 
and focusing on the simpler case of $\Z(\sigma )$.
(The same ideas apply for any generalized Zeta function $\Z(\sigma,v)$
and for complex $\sigma$, but the formulae get more involved.)
Here, the Riemann Hypothesis is {\sl de facto\/} implied throughout
(there being no numerical counter-example).

Numerically, $Z(\sigma)$ looks almost indistinguishable 
from $\Z(\sigma)$ for $\sigma \ge 0$,
because ${\tau_k}^2 \gg {1 \over 4} \ (\forall k)$
(an empirical fact; already, ${\tau_1}^2 \approx 199.790455$).
If we expand $Z(\sigma)$ in terms of $\Z(\sigma)$ according to
eq.(\ref{SHZ0}), and make the roughest approximations, we get that
$Z(\sigma) \approx \Z(\sigma) [1- \sigma / (4 {\tau_1}^2) ]$: i.e., 
the very first correction term is only of relative size $\approx \sigma/800$.
As other related numerical observations:
\begin{equation}
0 < \Z(\sigma) -Z(\sigma) < 0.0003 \mbox{ for all real }\sigma > 0 ;
\end{equation}
\begin{equation}
A \defi 4[Z'(0)-\Z'(0)] \approx -0.0231003495 \quad vs \quad 
B \equiv -Z(1)  \approx -0.0230957090 
\end{equation}
(not only is $A$ small, but moreover, 
$A \equiv 4 \log \Xi(\hf)$ by eq.(\ref{Z20}) 
and $B \equiv [\log \Xi]'(0)$ by eq.(\ref{Z11}),
hence $|B-A|<5 \times 10^{-6}$ reflects how {\sl little\/}
the function $\log \Xi(s)$ deviates from the parabolic shape $As(1-s)$ 
over the interval $[0,1]$). 

\medskip

We now focus on the numerical evaluation of $\Z(\sigma)$ itself.
The defining series (\ref{Z2}) converges more and more poorly 
as $ \sigma \to \hf^+$ (with divergence setting in at $\sigma =\hf$). 
We then replace a far tail of that series ($k>K$) by an integral 
according to the integrated-density estimate (\ref{GE}),
and formally obtain a kind of Euler--Maclaurin formula,
\begin{eqnarray}
\label{EMc}
\Z (\sigma ) \si=\sf \lim_{K \to +\infty} S_K (\sigma ), \qquad S_K (\sigma ) \defi
\sum_1^{K-1} {\tau_k}^{-2 \sigma }+ \hf {\tau_K}^{-2 \sigma } 
+ {\overline R}_K (\sigma ), \\
{\overline R}_K (\sigma ) \si\defi\sf \int_{\tau_K}^{+\infty} T^{-2 \sigma }
\d {\overline \mathcal N}(T) 
= {1 \over 2 \pi} {{\tau_K}^{1-2 \sigma } \over 2 \sigma -1} 
\Bigl[ \log {\tau_K \over 2 \pi}+{1 \over 2 \sigma -1}  \Bigr] \nonumber
\end{eqnarray}
(similar formulae can be written for $\Z'(\sigma)$, $\Z(\sigma,v)$,
finite parts at $\sigma=\hf$, etc.).

\begin{table}
\begin{tabular}  {ccc}
\hline \\[-12pt]
$\sigma$ & ~~ $Z(\sigma) = \sum\limits_{k=1}^\infty ({\tau_k}^2+\qt)^{-\sigma}
\quad [v=\qt]~~ $ &
~~ $\Z(\sigma) = \sum\limits_{k=1}^\infty {\tau_k}^{-2\sigma}
\quad [v=0] ~~ $ \\[10pt]
\hline \\[-10pt]
$-1$ & $-0.0625^\ast$ & $-0.28125^\ast$ \\
$-3/4$ & $1.69388$ & $0.54319$ \\
$- 1/4$ & $0.800805$ & $0.785321$  \\
$0$ & $0.875^\ast$ & $0.875^\ast$ \\
\hline \\[-10pt]
{\sl derivative at 0} & {\sl 0.8060429} & {\sl 0.8118179} \\[3pt]
$+1/4$ & $1.548829$ & $1.549060$ \\
{\sl finite part at 1/2} & {\sl 0.251546} & {\sl 0.251637} \\
$+3/4$ & $0.247730$ & $0.247760$ \\
$+1$ & $0.0230957$ &  $0.0231050$ \\
$+3/2$ &  $0.0007287$ &  $0.0007295$ \\
$+2$ & $0.0000371$ & $0.0000372$ \\
\hline
\end{tabular}
\caption{Numerical values for $\zeta$-Zeta functions of the Riemann zeros
($^\ast$: exact values).
Implied precision is expected to hold, but not guaranteed.
}
\end{table}

The approximate remainder term $\overline R_K(\sigma )$ balances 
the dominant {\sl trend\/} of the partial sums
$\sum\limits_{k<K}{\tau_k}^{-2\sigma }$.
It thereby accelerates the convergence of the partial sums in eq.(\ref{EMc}) 
for $\sigma > \hf$, 
while for $\sigma \le \hf$ it counters their divergent trend,
so that $S_K (\sigma )$ converges (as $K \to +\infty$) when $\sigma >0$ 
(\cite{G1}, p.116, last line).
The next obstruction to convergence arises at $\sigma=0$ but is of another type:
$S_K(\sigma ) $ displays {\sl erratic fluctuations\/} in $K$
({\sl roughly\/} of the order ${\tau_K}^{-2 \sigma } (\log \log \tau_K)^{1/2}$,
according to \cite{O}, eq.(2.5.7)),
and those numerically blow up indeed (as $K \to +\infty$) when $\sigma \le 0$.
Further convergence now requires to perform a damping of those fluctuations 
(as argued previously for ``chaotic" spectra \cite{BSV}).
Here, a {\sl Cesaro averaging\/} (defined by 
$\langle S \rangle_K \defi K^{-1} \sum\limits_1 ^K S_{K'} $) 
{\sl appears\/} to work well initially (results can be verified at $\sigma=0$), 
but not very far down:
already at $\sigma = -0.25$, the fluctuations of $\langle S \rangle_K(\sigma)$
itself retain a standard deviation $> 10^{-3}$ up to $K \approx 10^5$.
So, instead of pursuing ever more severe (and unproven, after all) 
numerical regularizations as $\sigma$ decreases below $\hf$,
we advocate the switch to the continuation formulae (\ref{ZC}--\ref{ZR1}) 
for numerical work as well. 
Thus, we first tested eq.(\ref{ZR1}) against eq.(\ref{EMc}) for $\Z(+\qt)$, 
then used it to evaluate $\Z(-\qt)$, 
plus eq.(\ref{SHZ0}) with $v=\qt$ (3 terms sufficed) to obtain $Z(-\qt)$.

Table 2 gives a summary of the numerical results we obtained.
(We found no earlier analogs, except for the {\sl other\/} special sums
${\mathscr Z}_n$ in \cite{Ms2}, \cite{L}.)

\end{document}